\documentclass[a4paper]{article}

\usepackage{amsfonts}
\usepackage{amssymb}
\usepackage{latexsym}

\newtheorem{teo}{Theorem}[section]

\newtheorem{lem}[teo]{Lemma}
\newtheorem{prop}[teo]{Proposition}
\newtheorem{cor}[teo]{Corollary}
\newtheorem{prob}[teo]{Problem}

\newtheorem{defin}[teo]{Definition}
\newcommand{\dok}[1]{{\bf Proof.} $\:$ {#1} $\qquad \square$ \medskip}

\begin{document}


\title{Varieties of representations of groups and varieties of
associative algebras}

\maketitle    

\begin{center}

\author{E. Aladova$^{a,}$\footnote{ The author is supported by the Israel
Science Foundation  grant number 1178/06. } ,
      \and B. Plotkin$^{b}$}

 \smallskip
        {\small
               $^{a}$ Department of Mathematics,
          Bar-Ilan University,

          52900, Ramat Gan, Israel

                {\it E-mail address:} aladovael@mail.ru
        }

\smallskip

 \smallskip
        {\small
               $^{b}$ Department of Mathematics,
          Hebrew University of Jerusalem,
          91904, Jerusalem, Israel

                {\it E-mail address:} plotkin@macs.biu.ac.il
        }
\end{center}





\begin{abstract}

This paper is tightly connected with the book~\cite{PV}. Since
this book had not been translated from Russian into English we
give here the brief review of the basic definitions and results
from~\cite{PV}.  We prove also
 new results of the same spirit.
 They are related to dimension subgroups,
 varieties of representations of groups and  varieties of associative algebras.
 The main emphasis is put on the varieties of representations of groups
  induced by the varieties of associative algebras. We provide the
reader with the list of some open problems.

\end{abstract}

\tableofcontents


\section{Introduction}

It is well known that there are close connections between
varieties of groups and varieties of representations of groups,
between varieties of algebras and varieties of representations of
algebras. The main goal of  this paper is to consider some relations
between varieties of representations of groups and varieties of
associative algebras. In particular, we consider polynomial
varieties of representations of groups which are related to varieties of associative algebras.

The set of
 all varieties of representations of groups is a free
semigroup under multiplication of varieties
\cite{Plotkin-Grinberg}. All polynomial varieties of
representations of groups form a subsemigroup  in this semigroup.
So we come up with the problem: is the semigroup of all polynomial
varieties of representations of groups free? To solve this problem
we are looking for some examples of polynomial and non-polynomial
varieties. For instance, the variety generated by the
representation $(K^{n},T_{n}(K))$ is a polynomial variety, while
the variety generated by the representation $(K^{n},UT_{n}(K)),$
is a non-polynomial variety. Here $K$ is a field, $T_{n}(K)$ is
the group of triangular matrices and $UT_{n}(K)$ is the group of
unitriangular matrices over $K$ which act naturally on the space
$K^{n}$.

The other our interest is concerned with dimension subgroup.  
There exist a lot of problems around classical dimension
subgroups (see \cite{Gupta}, \cite{MikhPassi},
\cite{MikhPassi_book}, \cite{Passi}, \cite{Rips} and references therein). We define a
concept of the dimension subgroup for varieties of representations of
groups. Here we consider the dimension subgroup of a group $G$ as
the intersection of kernels of all representations of the group $G$
in the given variety of representations of groups. In particular,
we especially distinguish the case of polynomial varieties of
representations of groups. We also give a concept of the dimension
subgroup for varieties of associative algebras, and consider some
problems concerning these subgroups.

A few words about the structure of the paper. In
Section~\ref{S1_VarOfReprGr} we give basic definitions: we define
the category {\bf Rep}-$K$ which is the category of all
representations of groups over a fixed ring $K$; the variety {\bf
Rep}-$K$, and consider some operations on the varieties, in
particular, triangular product of varieties of representations of
groups. In Section~\ref{S3_VarOfReprAlg} we recall  known
results about varieties of representations of algebras and
varieties of algebras from the book~\cite{PV}.
Some of these results which are very important for us we give with
the proofs. In Section~\ref{S4_VarOfReprGrVarAlg} we study
relations between varieties of representations of groups and
varieties of associative algebras. Here we deal with the following
problem: to find connection between interesting varieties of
representations of groups and some well-known varieties of
associative algebras. In the last section we give concepts of the
dimension subgroup for variety of representations of groups and
variety of associative algebras, and compare two approaches to the notion of the
dimension subgroup.

\section{Varieties of representations of groups. General
view}\label{S1_VarOfReprGr}

\subsection{Category and variety Rep-K}\label{SubS1.1_CatVarRep-K}


\paragraph{1.1.1 Category Rep-$K$. }

Let $K$ be a commutative and associative ring with unit. Consider
pairs $(V,G)$ where $V$ is a $K$-module and $G$ is a group acting
on the module $V$. In other words, for each $(V,G)$ there is a function from
$V\times G$ to $V$ defined by $(v,g)\to v\circ g$ satisfying the following
conditions:
\begin{enumerate}
\item
 $v_1\circ (g_1g_2)=(v_1\circ g_1)\circ g_2$;
\item
 $v_1\circ 1=v_1$, where 1 is the unit of the group $G$;
\item
 $(v_1+v_2)\circ g_1=v_1\circ g_1+v_2\circ g_1$;
\item
 $(\lambda v_1)\circ g_1=\lambda (v_1\circ g_1)$,
\end{enumerate}
where $v_1,v_2\in V$, $g_1,g_2\in G$, $\lambda \in K$.
 The pair $(V,G)$ can be viewed as two-sorted
algebra-representation over the ring $K$ (for terminology see books
\cite{Plotkin_Autom}, \cite{PV}).

On the other hand, if we have a group $G$ acting on a $K$-module
$V$ then there is a representation $\rho : G\to \mbox{Aut }V$ of
the group $G$ in the $K$-module $V$ defined as follows $g^{\rho}(v)=v\circ
g$, where $g\in G$, $v\in V$, $g^{\rho}$ is the image of $g$ under
the homomorphism $\rho$.

Let $(V,G)$ be a  representation. Then $V$ can be viewed as a
$G$-module. Moreover, since the group $G$ acts on the $K$-module
$V$ then the group algebra $KG$ also acts on $V$ and,
consequently, $V$ is a $KG$-module.

The class of all pairs $(V,G)$ over a fixed ring $K$ forms the
category {\bf Rep}-$K$. Objects of this category are
representations $(V,G)$, morphisms are homomorphisms $\mu=(\alpha,
\beta )$ of representations defined as follows:

\begin{defin}
A homomorphism of a representation $(V_1, G_1)$ to a
representation $(V_2,G_2)$ is a pair of homomorphisms
$$
 \mu=(\alpha, \beta ): (V_1,G_1)\to (V_2,G_2),
$$
 where $\alpha :
V_1\to V_2$ is a homomorphism of  $K$-modules, $\beta : G_1\to
G_2$ is a homomorphism of  groups and the following condition
holds
$$
(v\circ g)^{\alpha}=v^{\alpha}\circ g^{\beta}.
$$
\end{defin}

\begin{defin}
A representation $(U, H)$ is a subrepresentation  of a
representation $(V, G)$ if $H$ is a subgroup of the group $G$, $U$
is a $K$-submodule of $V$ and the restriction of the action $G$ on
$V$ coincides with the action $H$ on $U$.
\end{defin}

Along with the subrepresentations we consider a concept of
over\-representation. The defined  below induced representations is a particular case of  this concept.

 Let $H$ be a subgroup of a group $G$, and
let $(U,H)$ be a representation of the group $H$. Consider the
following category: objects of this category are $H$-homomorphisms
$\varphi : U\to V$ of the $H$-module $U$ to 
$G$-modules $V$, morphisms of an object $\varphi_1 : U\to V_1$
to an object $\varphi_2: U\to V_2$ are $G$-homomorphisms $\psi :
V_1\to V_2$ subject to the  following commutative diagram

\begin{picture}(100,60)\
\put(55, 20){\vector(2,1){40}} \put(40,20){$U$} \put(100,
35){$V_1$} \put(70,35){$\varphi_1 $}%
\put(55, 20){\vector(2,-1){40}} \put(100,-5){$V_2$} \put(70,-5){$\varphi_2 $}%
\put(105, 30){\vector(0,-1){20}} \put(110, 15){$\psi$}
\end{picture}

\bigskip

 Let $\tau : U\to V$ be an initial object of this category.
A representation $(V,G)$ is called  induced
by the representation $(U,H)$ with the given embedding of the
group $H$ to the group $G$. Note that $U$ is a right $KH$-module
and $KG$ is a left $KG$-module. Then we can consider the module
$V=U\otimes_{KH} KG$ and the representation $(V,G)$ is the
representation induced by the representation $(U,H)$.

There are products and coproducts of objects in the category {\bf
Rep}-$K$~\cite{MacLane}.

\begin{defin}
A representation $(V,G)$ is called the Cartesian product of
representations $(V_i,G_i)$ if $V$ is a Cartesian product of the
$K$-modules $V_i$, $G$ is a Cartesian product of the groups $G_i$
and for all $v\in V$ and all $g\in G$ we have $(v\circ
g)(i)=v(i)\circ g(i)$.
\end{defin}

The Cartesian product $(V,G)$ is a product of the objects
$(V_i,G_i)$ in the  category {\bf Rep}-$K$.

\begin{defin}
A representation $(V,G)$ is called the free
 product of the representations $(V_i,G_i)$ if $G$ is the free product of the groups $G_i$,
 $V$ is the direct sum of $G$-submodules $U_i$ such that representations $(U_i,G)$ are induced
 by the representations $(V_i,G_i)$ with the
natural embedding $G_i\to G$.
\end{defin}

The free product $(V,G)$ is
 a coproduct of the objects $(V_i,G_i)$ in the category {\bf
Rep}-$K$.

Now we describe a free objects in the category {\bf Rep}-$K$. Let
$(V,G)$ be an arbitrary representation. Let $Y$ be a subset of $G$
and let $X$ be a subset of $V$. We say that the pair of the sets
$\{X,Y\}$ generates the representation $(V,G)$ if the group $G$ is
generated by the set $Y$ and the module $V$ is generated by the
set $X$ as $KG$-module.

\begin{defin}
A representation $(W,F)$ is called a free representation over a
pair of the sets $\{X,Y\}$ if for every representation $(V,G)$ any
pair of the functions $\mu_1:Y\to G$ and $\mu_2:X\to V$
can be extended  uniquely up to the homomorphism of the
representations $(W,F)\to (V,G)$.
\end{defin}

\paragraph{1.1.2 Variety Rep-$K$. }

Along with the category {\bf Rep}-$K$ we  consider variety of
representations of groups {\bf Rep}-$K$. This variety consists of
all two-sorted algebras-representations. Identities of the variety
{\bf Rep}-$K$ are identities of groups, identities of $K$-modules
and identities of action written earlier.

Describe the free representations in the  variety {\bf Rep}-$K$. Let
$\{ X, Y\}$ be a pair of sets.  Let $F=F(Y)$ be the free group over
the set $Y$, $W=W(X,Y)$ be the free $KF$-module over the set $X$.
Moreover, we have
$$
 W(X,Y)=X\cdot KF=\bigoplus\limits_{x\in X} x \cdot KF,
$$
and each element $w\in W$ have the form:
$$
w=x_1u_1+\dots +x_nu_n,
$$
where  $u_i=u_i(y_1,\dots ,y_m)\in KF$, $x_i\in X$. Define an
action of the group $F$ on the module $W$ as follows:
$$
w\circ f=(x_1u_1+\dots +x_nu_n)\circ f= x_1(u_1f)+\dots
+x_n(u_nf),
$$
where $w\in W,\ f\in F,\ u_if\in KF$.

\begin{defin}
A representation $(W,F)$ with the action $w\circ f$ above  is the
free representation over $\{ X, Y\}$.
\end{defin}

So, we can speak about two types of identities of representations:
$$
 w\equiv 0 \mbox{ is an identity of  action},\  w\in W(X,Y),
$$
 and
$$
 f\equiv 1 \mbox{ is a group identity,} \ f\in F.
$$
 In fact we will consider only identities of the action.
 Moreover, it is sufficient to consider identities of the form
$$
x\circ u\equiv 0,\ u\in KF,\  x\in X.
$$

Let $F$ be a free group of countable rank with free generators
$y_1,y_2,\dots $, and let $KF$ be a group algebra of $F$. Let
$u(y_1,\dots ,y_n)$ be an element of $KF$.

\begin{defin}
We say that a representation $(V,G)$ satisfies an identity $x\circ
u(y_1,\dots ,y_n)\equiv 0$  if for all $v\in V$ and all $g_i\in G$
we have $v\circ u(g_1,\dots ,g_n)=0$.
\end{defin}


Consider some examples of identities of representations of groups.

\bigskip
 {\bf Examples.}
\begin{enumerate}
\item
 A representation $(V,G)$ satisfies the identity
 $$
  x\circ (y-1)\equiv 0
 $$
 if the group $G$ acts trivially on the module $V$.

\item
 A representation $(V,G)$ satisfies the identity
 $$
  x\circ (y_1-1)(y_2-1)\dots (y_n-1)\equiv 0
 $$
 if the module $V$ has a $G$-invariant series of
 submodules of the length $n$:
 $$
  0=V_0\subset V_1\subset \dots \subset V_n=V,
 $$
 such that the
 group $G$ acts trivially on the factors. The representation
 $(K^{n},UT_{n}(K))$ satisfies this identity. Here $UT_{n}(K)$ is
 the group of unitriangular matrices over a field $K$, and this group acts in a natural
 way on the module $K^{n}$.

\item
 Let $(V,G)$ be a finite dimension representation over a field
 $K$ and  let $n=\mbox{dim}\ V$. It is known that the algebra
 $\mbox{End}V$ satisfies the standard identity of Amitzur-Levitzky of the degree $2n$, that is
 $$
  s_{2n}=\sum_{\sigma\in S_{2n}}
  (-1)^{\mbox{sgn}\sigma}x_{\sigma(1)}\dots x_{\sigma(2n)}\equiv 0,
 $$
 where $S_{2n}$ is the symmetric group. So all
 finite dimension representations possess non-trivial identities.


\end{enumerate}

\subsection{Varieties of representations as subvarieties of the Rep-K}\label{SubS1.2_VarSubvarRep-K}

In this part of the paper we will speak about varieties of
representations which are subvarieties of {\bf Rep}-$K$. This means that varieties in question
satisfy a system of identities of the form $x\circ u\equiv 0$. We
suppose that the free group $F=F(Y)$ is a group 
of
countable rank.

Let $\mathfrak X$ be a class of representations from the variety
{\bf Rep}-$K$. Let $I_{\mathfrak X}(KF)$ be a set of all $u\in KF$
such that the identity $x\circ u\equiv 0$ holds in every
representation from $\mathfrak X$. It easy to see that the set
$I_{\mathfrak X}(KF)$ is a two-sided ideal of the group algebra
$KF$ which is invariant under all endomorphisms of the group $F$, i.e.,
 $I_{\mathfrak X}(KF)$ is a fully invariant ideal of
$KF$. Furthermore there is a one-to-one correspondence between
 the varieties of representations of groups over $K$ and
 the fully invariant ideals of $KF$~\cite{PV}. 

Let $\mathfrak X$ be a class of representations of groups over
$K$. Denote by $\mbox{Var}\mathfrak X$ the variety of
representations generated by the class $\mathfrak X$.
For varieties of representations of groups we have the following
invariant description~\cite{PV}:
\begin{teo}
A class of representations of groups forms a variety of
representations of groups if and only it is closed under
subrepresentations, homomorphic images, Cartesian
products and 
saturations.
\end{teo}
Remind that a class of representations of groups $\mathfrak X$ is
closed under the saturation if the following condition is true: if
a representation $(V,H)$ lies in the class $\mathfrak X$ then all
representations $(V,G)$ such that $(V,H)$ is a right epimorphic
image of the representation $(V,G)$ also belong to $\mathfrak X$.


\subsection{Algebra of varieties}\label{SubS1.3_AlgOfVar}

\paragraph{1.3.1 Operations on varieties.}

Now we consider the algebra of subvarieties of {\bf Rep}-$K$.

Let $\mathfrak X_1$ and $\mathfrak X_2$ be two subvarieties from
{\bf Rep}-$K$.

\begin{defin}[\cite{PV}]
A product of  varieties  $\mathfrak X_1$ and $\mathfrak X_2$ is a
variety $\mathfrak X_1\mathfrak X_2$ defined as follows: a
representation $(V,G)$ belongs to $\mathfrak X_1\mathfrak X_2$ if
and only if the module $V$ has a $G$-submodule $U$ such that
$(U,G)\in \mathfrak X_1$ and $(V/U,G)\in \mathfrak X_2$.
\end{defin}

The operation of the multiplication of varieties is associative.
Thus we have a semigroup $\mathfrak M=\mathfrak M(K)$ of varieties
of representations of groups.

For example,
%
 let $\mathfrak S$ be a variety satisfying the identity $x\circ
 (y-1)\equiv 0$ and $\mathfrak S^{n}$ be a variety determined by the
 identity $x\circ (y_1-1)\dots (y_n-1)\equiv 0$. Then the variety
 $\mathfrak S^{n}$ is 
 $n$-th degree of the variety $\mathfrak S$.

\bigskip
There is the following theorem~\cite{Plotkin-Grinberg} (see
also~\cite{PV}):

\begin{teo}\label{Th_1_1}
If $K$ is a field then the semigroup $\mathfrak M=\mathfrak M(K)$
is a free semigroup. $\qquad \square$
\end{teo}

In the proof of Theorem \ref{Th_1_1} the construction of
triangular product of representations is used. Define this
structure (see \cite{PV}, \cite{Vovsi}).
 Let $(V_1,G_1)$ and $(V_2,G_2)$  be representations. Define a representation
 $(V,G)$ as follows
$$
 (V,G)=(V_1,G_1)\bigtriangledown (V_2,G_2),
$$
where $V=V_1\oplus V_2$ is the direct sum of the modules
and elements of the acting group $G$ are matrices of the form:
$$
g=
\left(
\begin{array}{cc}
g_1 & \varphi \\
0   & g_2     \\
\end{array}
\right),
$$
where $g_1\in G_1$, $g_2\in G_2$, $\varphi \in
\mbox{Hom}(V_1,V_2)$.

Let
$$
g= \left(
\begin{array}{cc}
g_1 & \varphi \\
0   & g_2     \\
\end{array}
\right),\ \
 g'= \left(
\begin{array}{cc}
g'_1 & \varphi ' \\
0   & g'_2     \\
\end{array}
\right)
$$
be elements of the group $G$. The product of these elements is defined as follows:
$$
gg'= \left(
\begin{array}{cc}
g_1g'_2 & \overline{g_1}\varphi '+\varphi \overline{g'_2} \\
0   & g_2g'_2     \\
\end{array}
\right),
$$
where $g_i, g'_i\in G_i$, $\varphi, \varphi '\in
\mbox{Hom}(V_1,V_2)$, $\overline{g_1}\in \overline G_1$,
$\overline{g'_2}\in \overline G_2$ and $(V_i,\overline G_i)$ are
faithful representations corresponding to the representations
$(V_i, G_i)$, $i=1,2$.

The group $G$ acts on the module $V$ by the rule:
$$
(v_1+v_2)\circ
 \left(
 \begin{array}{cc}
 g_1 & \varphi  \\
 0   & g_2     \\
 \end{array}
 \right)=
 v_1\circ g_1+\varphi (v_1)\circ g_2+v_2\circ g_2,
$$
where $v_1\in V_1$, $v_2\in V_2$.

Note that the submodule $V_1$ of the module $V$ is invariant under
the action of the group $G$. Moreover the representation $(V_1,G)$
has the same identities as the representation $(V_1,G_1)$ and the
representation $(V/V_1,G)$ has the same identities as the
representation $(V_2,G_2)$.

Henceforward $K$ is a field. The basic role in the proof of Theorem \ref{Th_1_1}
plays the following proposition.

\begin{prop}[\cite{PV}]\label{Prop_1-1} 
Over a field we have the following property:
$$
\mbox{Var }((V_1,G_1)\bigtriangledown (V_2,G_2))=\mbox{Var
}(V_1,G_1)\cdot \mbox{Var }(V_2,G_2). \  \qquad \square
$$
\end{prop}

Let $\Theta_1$ and $\Theta_2$ be varieties of groups.

\begin{defin} 
A product of varieties of groups $\Theta_1$ and $\Theta_2$ is a
class $\Theta_1\Theta_2$ defined as follows: a group $G$ belongs to
the variety $\Theta_1\Theta_2$ if it has an invariant subgroup $H$
such that $H\in \Theta_1$ and $G/H\in \Theta_2$.
\end{defin}

Note that the class $\Theta_1\Theta_2$ is a variety of groups. Let
$\mathfrak M^{g}$ be the semigroup of varieties of groups. It is
well-known that $\mathfrak M^{g}$ is a free semigroup
\cite{Neumann_3}, \cite{Shmel'kin}.

Now we assign an action of the semigroup $\mathfrak M^{g}$ on the
semigroup $\mathfrak M$ of varieties of representations of groups.

Let $\mathfrak X$ be a variety of representations of groups and
$\Theta $ be a variety of groups. Define a class of
representations of groups $\mathfrak X\times \Theta$  as follows:

\begin{defin}[\cite{PV}] 
A representation $(V,G)$ belongs to $\mathfrak X\times \Theta$ if
the group $G$ has a normal subgroup $H$ such that $(V,H)\in
\mathfrak X$ and $G/H\in \Theta$.
\end{defin}

Note that the class $\mathfrak X\times \Theta$ forms a variety of
representations of groups and there is the following

\begin{prop}[\cite{PV}]\label{Prop_1-2} 
Let $\mathfrak X_1$ and $\mathfrak X_2$ be varieties of
representations of groups and let $\Theta_1 $ and $\Theta_2$ be
varieties of groups. Then
\begin{enumerate}
\item
 $\mathfrak X_1\times \Theta_1\Theta_2=(\mathfrak
 X_1\times\Theta_1)\times\Theta_2$,

\item
 $\mathfrak X_1\mathfrak X_2\times \Theta_1=(\mathfrak
 X_1\times\Theta_1)(\mathfrak X_2\times\Theta_1)$. $\qquad \square$

\end{enumerate}
\end{prop}

The following theorem is a main theorem of this theory.

\begin{teo}[\cite{Plotkin_VarAndVarOfPairs}]\label{Th_1-2} 
If $K$ is a field then the semigroup $\mathfrak M^{g}$ acts freely
on the semigroup $\mathfrak M=\mathfrak M(K)$. $\qquad \square$
\end{teo}

In particular, this theorem yields  the following property: for
every variety of representations of groups $\mathfrak X$ there is
unique decomposition
$$
 \mathfrak X=(\mathfrak X_1\times \Theta_1)\dots (\mathfrak X_n\times
 \Theta_n),
$$
where all $\mathfrak X_i$ are indecomposable as $(\mathfrak X_{i_1}\times \Theta_{i_1})$. 

\bigskip

{\bf Example.}
 Let $\mathfrak S_n$ be the variety determined by the identity $x\circ (y-1)^{n}\equiv
 0$. The variety $\mathfrak S_n$ is indecomposable.

\bigskip

In the proof of Theorem~\ref{Th_1-2}  the construction of 
the wreath product of a representation of a group and a group is used.

Let $(V,H)$ be a representation  and $G$ be a group. Define 
$$
 (V,H)\mbox{ wr }G=(V^{G},H\mbox{ wr }G)=(V^{G},H^{G}\leftthreetimes G).
$$

Only in special cases we have the following property:
$$
 \mbox{Var }((V,G)\mbox{ wr }G)=\mbox{Var }(V,G)\times \mbox{Var
 }G.
$$

\paragraph{1.3.2 Varieties of groups and varieties of representations of groups.}

Now we give some remarks about connections between varieties of
groups and varieties of representations of groups.

Let $\Theta$ be a variety of groups. Define a class of
representations of groups $\omega\Theta$ as follows

\begin{defin} 
Let $(V,G)$ be a representation and $(V,\overline {G})$ the corresponding faithful representation. A representation $(V,G)$ lies in $\omega\Theta$ if the group $\overline {G}$
belongs to $\Theta$.
\end{defin}

Note that class $\omega\Theta$ is a variety of representations
of groups. It is obvious that $\omega\Theta=\mathfrak S\times
\Theta$, where $\mathfrak S$ is the variety of representations satisfying the
identity $x\circ (y-1)\equiv 0$.

Let $\mathfrak A$ be a variety of abelian group, $\mathfrak A_m$
be a variety  of abelian groups of exponent $m$. Let $K$ be an
integral domain and let $K^{*}$ be a multiplicative group of $K$.
Then we have two propositions (see~\cite{PV}):

\begin{prop}  
If $K^{*}$ is infinite then the variety of representations of
groups $\omega\mathfrak
 A$ is generated by the representation $(K,K^{*})$, where $K^{*}$
 acts on the cyclic module $K$ by multiplication.
That is
$$
\omega\mathfrak A=\mbox{Var}(K,K^{*}).
$$
\end{prop}

\begin{prop}  
If the order of $K^{*}$ is equal $m$ then the variety of
representations of groups $\omega\mathfrak
 A_m$ is generated by the representation $(K,K^{*})$.
\end{prop}

 Using the proposition~\ref{Prop_1-2} we have the following fact: 
 $$
  (\omega\mathfrak A)^{n}=(\mathfrak S\times \mathfrak A)^{n}=\mathfrak S^{n}\times
  \mathfrak A.
 $$
 The variety $\mathfrak S^{n}$ is generated by the representation
 $(K^{n},UT_{n}(K))$ (see ~\cite{Grinberg}),
 and the variety $\mathfrak S^{n}\times \mathfrak A $ is generated by the representation
 $(K^{n},T_{n}(K))$~(\cite{Krop}, see also \cite{PV}).

 Now define the operator $\omega '$ from
 varieties of groups $\Theta$ to varieties of representations of groups  $\omega '\Theta$.
Let $(V,G)$ be a representation of a group. Consider $V$ as an abelian group and take  the
 semidirect product of groups
$V\leftthreetimes G$. 
 We say that the representation $(V,G)$ belongs to the class
$\omega '\Theta$ if $V\leftthreetimes G\in \Theta$. 
The
class $\omega '\Theta$ is a variety of representations of
groups.
 Here the symbol $'$ is interrelated with the following fact: the
identities of the variety $\omega '\Theta$ are given by the special
Fox derivatives used in the theory of formal power series \cite{Fox}. We
have the following conditions:
\begin{enumerate}
\item
 $\omega '(\Theta_1\Theta_2)=((\omega '\Theta_1)\times \Theta_2)\cdot \omega
 '\Theta_2$, 
\item
 $\omega '[\Theta_1,\Theta_2]=(\omega \Theta_1\cdot \omega '\Theta_2)\cap (\omega\Theta_2\cdot \omega
 '\Theta_1)$, 
\end{enumerate}
\noindent
where $[\Theta_1,\Theta_2]$ is the commutator of varieties.
Second formula confirms the "differential character" of the
operator $\omega '$.

In the sequel we will give some other  remarks and
propositions concerning this theory.

\section{Varieties of representations of associative algebras}\label{S3_VarOfReprAlg}

In this section we consider some definitions and facts concerning
varieties of representations of associative algebras. Here we rely
on
 $\S$14 of the book~\cite{PV}
since
there is an essential connections of these facts with  
new results.

Let $K$ 
be an associative and commutative ring with unit, let $A$ be an
associative algebra over $K$ and let $V$ be a $K$- and $A$-module.
Let $(V,A)$ be a representation of the algebra $A$. We also can
consider the category of representations of associative algebras
over a fixed ring $K$ and various algebras $A$. A variety $\Sigma$
in this category can be defined as a class of representations of
algebras closed under the subrepresentations, homomorphic images,
Cartesian products and saturations. These operators are defined
in the same way as for representations of groups (see section
\ref{SubS1.1_CatVarRep-K}, \ref{SubS1.2_VarSubvarRep-K}).
Similarly we can define the product of two representations of
algebras over a fixed ring $K$ (see section
\ref{SubS1.3_AlgOfVar}) and so we have the semigroup of
representations of algebras $\mathfrak M^{a}=\mathfrak M^{a}(K)$.

Let $\Sigma$ be a variety of  associative $K$-algebras. Define a
class of representations of algebras $\omega\Sigma$ as follows:

\begin{defin}
A representation $(V,A)$ belongs to $\omega\Sigma$ if the algebra
$\overline A, $ where $(V,\overline A)$ is faithful 
belongs to
$\Sigma$.
\end{defin}


 The class of representations $\omega\Sigma$ is
a variety of representations. A variety  of representations of algebras of the form  $(0,A)$, where $0$ is the zero algebra is called a singular variety.
Let $\mathfrak Y$ be a non-singular variety  of representations of
algebras.


\begin{defin}\label{def_FaithulRep}
 Denote by   $\overrightarrow{\mathfrak Y}$ the class of all algebras
 admitting a faithful representation in the variety $\mathfrak Y$.
\end{defin}

We will prove that the class $\overrightarrow{\mathfrak
Y}$ is a variety of algebras. Note that this fact makes difference between the cases of algebras and groups.
It is known
\cite{PV} that a class $\overrightarrow{\mathfrak X}$ of all
groups admitting a faithful representation in the variety of
representations of groups $\mathfrak X$ constitutes a quasivariety of
groups.

 Let $A$ be a $K$-algebra with unit and $J$ be
an ideal of $K$ 
 annihilating the $K$-module $A$. We
can consider the algebra $A$ as a $K/J$-algebra.

Let $\mathfrak Y$ be a variety of representations of algebras over
a ring $K$ and let $\mathbb F$ be a free $K$-algebra with unit.
Let $I_{\mathfrak Y}(\mathbb F)\triangleleft \mathbb F$ be the
ideal of identities of the variety $\mathfrak Y$. So the
representation $(\mathbb F/I_{\mathfrak Y}(\mathbb F),\mathbb F)$
is the free representation of the variety $\mathfrak Y$. Let $J$ be
an ideal of $K$ 
of all $\alpha\in K$ such that
$\alpha$ annihilates $K$-module $\mathbb F/I_{\mathfrak Y}(\mathbb
F)$, i.e.,  $\alpha\in J$ if $\alpha\cdot 1\in I_{\mathfrak
Y}(\mathbb F)$. It is clear that, if $(V, A)$ is a representation from
the variety $\mathfrak Y$ then the ideal $J$ annihilates the
module $V$. Moreover, if $A$ is an algebra from
$\overrightarrow{\mathfrak Y}$ then the ideal $J$ annihilates $A$.
The ideal $J$ is called a $K$-annihilator of the variety
$\mathfrak Y$.

We have the following

\begin{lem}[\cite{PV}]\label{Lem_3-1}
 Let $\mathfrak Y$ be a variety of representations of algebras. Then for every
 algebra $A$ from the variety $\overrightarrow{\mathfrak Y}$ the regular
 representation $(A,A)$ belongs to the variety $\mathfrak
 Y$.
\end{lem}

\dok{
 Since $A\in \overrightarrow{\mathfrak Y}$ then there exists a faithful
 representation $(V,A)$ from $\mathfrak Y$. Let $U=K/J$,
 where $J$ is the $K$-annihilator of $\mathfrak Y$. The
 representation $(U,0)$ also belongs to $\mathfrak Y$ as a
 subrepresentation of the free representation $(\mathbb F/I_{\mathfrak Y}(\mathbb F),\mathbb
 F)$. Using the homomorphism $A\to 0$ we define a zero
 representation $(U,A)$ and $(U,A)\in \mathfrak Y$.

 Let us consider the $A$-module $V\oplus U$.
 The representation $(V\oplus U,A)$ is  a faithful
 representation and it belongs to $\mathfrak Y$.






 Now we consider the representation $(A,A)$.
 Let $w$ be an element of $V\oplus U$. Define the map $\mu_{w}: A\to V\oplus
 U$ by the rule $\mu_{w}(a)=w\circ a$, where $a\in A$.
 The map $\mu_{w}$ is a homomorphism of $K$- and $A$-modules.
 Let $U_{w}=\mbox{Ker }\mu_{w}$. Then $A/U_{w}\le V\oplus
 U$ and the representation $(A/U_{w}, A)$
 belongs to the variety $\mathfrak Y$ as a subrepresentation  of
 the representation $(V\oplus U, A)\in \mathfrak Y$.
 Since $\mathfrak Y$ is a variety then the representation $(\prod_{w}(A/B_{w}), A)$
 belongs to the variety $\mathfrak Y$.
 The algebra
 $A$ acts faithfully on $V\oplus U$ and $\bigcap
 U_{w}=\{0\}$. Using the Remak's theorem we have $A\le \prod_{w}(A/U_{w})$.
 Thus the representation $(A,A)$ belongs to $\mathfrak Y$ as a subrepresentation of
 $(\prod_{w}(A/U_{w}),A)$.
} 

We have the following

\begin{prop}[\cite{PV}]\label{Prop_3-1} 
Let $\mathfrak Y$ be a variety of representations of associative
$K$-algebras. The class $\overrightarrow{\mathfrak Y}$ is a
variety of associative algebras.
\end{prop}

\dok{
 It is necessary to prove that the class $\mathfrak Y$ closed under
 subalgebras, Cartesian products and homomorphic images.

 It is obvious that first two conditions are hold.
 We show that the class $\overrightarrow{\mathfrak Y}$ closed under homomorphic images.
 Let $A$ be an algebra from the class $\overrightarrow{\mathfrak Y}$ and
 let $\mu: A\to A_1$ be an epimorphism of algebras.
 The
 epimorphism $\mu$ induces an epimorphism of the representations
 $\overline\mu : (A, A)\to (A_1, A_1)$.
 Since the regular representation $(A,A)$ belongs to $\mathfrak
 Y$ then the regular representation $(A_1,A_1)$ also belongs to $\mathfrak
 Y$. This representation is faithful,
%
%
 so the algebra $A_1$ belongs to $\overrightarrow{\mathfrak
 Y}$. Thus the class $\overrightarrow{\mathfrak Y}$ forms a variety of algebras.
} 

Thus the map $\overrightarrow{}\ $
 is a map from the set of all varieties of
representations of algebras over $K$ to the set of all varieties
of algebras. This fact is not true for varieties of
representations of groups and for varieties of representations of
Lie algebras (see~\cite{PV}). But we can note that always there is
the following property
$\overrightarrow{\omega\Sigma}=\Sigma$ (where $\Sigma$ is a
variety of groups, algebras or Lie algebras).

Let $\Sigma$ be a variety of associative algebras over a ring $K$,
let $J$ be an ideal of $K$ annihilating all algebras from the
variety $\Sigma$. Denote by $\omega_{J}\Sigma$ an intersection of
the variety of representations  $\omega\Sigma$ and the variety of
all representations $(V,A)$ such that the ideal $J$ annihilates
the module $V$. Then the following proposition holds:

\begin{prop}[\cite{PV}]\label{Prop_3-2} 
Let $\mathfrak Y$ be a variety of representation of algebras over
a ring $K$, let $J$ be the $K$-annihilator of the variety
$\mathfrak Y$. Then
$$
 \mathfrak Y=\omega_{J}\overrightarrow{\mathfrak Y}.
$$
\end{prop}

\dok{
 It is not difficult to see that $\mathfrak Y\subseteq \omega_{J}\overrightarrow{\mathfrak
 Y}$. We prove that $\omega_{J}\overrightarrow{\mathfrak Y}\subseteq \mathfrak Y$.
 Let $(V,A)$ be a representation from $\omega_{J}\overrightarrow{\mathfrak Y}$. We consider two
 cases: $V$ is a cyclic $A$-module and general case.

 Let $V$ be a cyclic $A$-module and let $(V,\overline A)$
 be the faithful representation corresponding the representation
 $(V,A)$. Then $\overline A\in \overrightarrow{\mathfrak Y}$.
 According
 to the lemma~\ref{Lem_3-1} the regular representation
 $(\overline A,\overline A)$ belongs to $\mathfrak
 Y$. The representation $(V,\overline A)$ is a homomorphic image of
 the representation $(\overline A,\overline A)$,
 and then the representation $(V,\overline A)$ belongs
 to the variety $\mathfrak Y$.
 Since  $(V,\overline A)$ is a right epimorphic image of the representation $(V, A)$
 and $(V,\overline A)\in\mathfrak Y$
 then  $(V, A)$ belongs to $\mathfrak Y$.

 Now let $(V, A)$ be an arbitrary representation from $\omega_{J}\overrightarrow{\mathfrak
 Y}$. Since the $A$-module $V$ is generated by its cyclic
 submodules $V_i$ and the representations $(V_i, A)$
 belongs to $\mathfrak Y$
  then the representation $(V,A)$ belongs to $\mathfrak Y$.
} 

We have the following

\begin{teo}[\cite{PV}]\label{Th_3-1} 
Let $K$ be a field. Then there is a one-to-one correspondence
between varieties of $K$-algebras and varieties of representations
of algebras over $K$.
\end{teo}

\dok{
 From the propositions~\ref{Prop_3-1} and \ref{Prop_3-2}  follows
 that if $K$ is a field then the maps $\omega$ and $\overrightarrow{}\ $
 are mutually inverse. So the statement of the theorem is true.
} 

Now we give some remarks about identities of varieties of
representations of algebras.

Let $(V, A)$ be a faithful representation over a ring $K$ from a
variety $\mathfrak Y$, and let $\mathbb F$ be the free
$K$-algebra. An element $u\in\mathbb F$ is an identity of the
faithful representation $(V,A)$ if and only if $u$ is the identity
of the algebra $A$.

Thus, if $\mathfrak Y$ is a variety of representations of algebras
over a ring $K$ and $\mathfrak Y=\omega_{J}\Sigma$, where $\Sigma
=\overrightarrow{\mathfrak Y}$, then the ideal of identities of
the variety $\mathfrak Y$ is generated by the ideal of identities
of the variety $\Sigma$ and the elements $\lambda \cdot 1$, where
$\lambda \in J$.








Note that, if $K$ is a field then the ideal of identities of a
variety of representations of algebras $\mathfrak Y$ coincides
with the ideal of identities of the variety of algebras
$\Sigma=\overrightarrow{\mathfrak Y}$. 

Now let $K$ be a field, let $\mathbb F$ be a free associative
$K$-algebra of countable rank and let $\Sigma$ be a class of
associative $K$-algebras. Take  $u_i\in \mathbb F$  such
that $u_i\equiv 0$ are identities in all algebras from the class
$\Sigma$. Then the elements $u_i$ form an ideal of the algebra
$\mathbb F$ closed under all endomorphisms of $\mathbb F$. This
ideal is called a fully invariant ideal (or T-ideal) and moreover
there is  one-to-one correspondence between fully invariant ideals
of the algebra $\mathbb F$ and varieties of associative algebras
over a field $K$.

Let $\Sigma_1$ and $\Sigma_2$ be varieties of algebras over a
field $K$, let $I_1$ and $I_2$ be fully invariant ideals (T-ideals)
of the varieties $\Sigma_1$ and $\Sigma_2$ accordingly. Let
consider varieties of representations of algebras $\omega\Sigma_1$
and $\omega\Sigma_2$ and their product
$\omega\Sigma_1\cdot\omega\Sigma_2$. Let
$\Sigma=\omega^{-1}(\omega\Sigma_1\cdot\omega\Sigma_2)$ and let
$I$ be an ideal of identities of the variety $\Sigma$. It is not
difficult to check that $I=I_2\cdot I_1$.











So the semigroup $\mathfrak M^{a}(K)$ of all varieties of
representations of algebras over a given field $K$ is
anti-isomorphic to the semigroup of all $T$-ideals of the free
$K$-algebra $\mathbb F$. The well-known theorem of Bergman and Lewin
\cite{Bergman-Lewin} states that the semigroup of all $T$-ideals
of $\mathbb F$ is a free semigroup. Note, that this theorem can be
proved with a help of the construction of the triangular product
(see section \ref{SubS1.3_AlgOfVar}).

\section{Varieties of representation of groups associated with
varieties of associative algebras}\label{S4_VarOfReprGrVarAlg}

Let $K$ be a field and let $\Sigma$ be a variety of associative
$K$-algebras with unit. Define a class $\eta\Sigma$ of
representations of groups as follows:

\begin{defin}
A representation $(V,G)$ belongs to $\eta\Sigma$ if the
representation $(V,KG)$ lies in the variety $\omega\Sigma$.
\end{defin}

In other words, if $(V,\overline G)$ is the faithful
representation of group $\overline G$ corresponding to the
representation $(V,G)$ then the linear span of $\overline G$ in
the algebra $\mbox{End}V$ is an algebra from the variety $\Sigma$.

We have the following

\begin{prop}[\cite{PV}]\label{Prop_4_1} 
The class $\eta\Sigma$ is a variety of representations of groups.
\end{prop}

\dok{
 Let $(V,G)$ be a representation from the variety $\eta\Sigma$ and
 let $(U,H)$ be a subrepresentation of $(V,G)$. Then $(U,KH)$ is
 a subrepresentation of $(V,KG)$. The representation $(V,KG)$
 belongs to variety $\omega\Sigma$, so $(U,KH)$ also belongs to
 $\omega\Sigma$. Hence $(U,H)$ is a representation from the
 variety $\eta\Sigma$. Thus the class $\eta\Sigma$ is closed
 under subrepresentations. In a similar way it can be proved
 that the class $\eta\Sigma$ is closed under homomorphic
 images.

 Now we prove that the class $\eta\Sigma$ is closed under saturations.
 Let $(V,G)$ be a representation such
 that the faithful representation $(V,\overline G)$ belongs to
 $\eta\Sigma$. Then $(V,K\overline G)$ lies in $\omega\Sigma$ and
 the faithful representations for $(V,KG)$ and $(V,K\overline G)$ are
 isomorphic. So $(V,KG)$ belongs to $\omega\Sigma$ and the representation $(V,G)$ belongs to
 $\eta\Sigma$.

 It remains to prove that the class $\eta\Sigma$ is closed under
 Cartesian products. Let $(V,G)$ be a Cartesian product of the
 representations $(V_i,G_i)$, $i\in I$, and $(V_i,G_i)\in
 \eta\Sigma$. Let $(V,A)$ be a Cartesian product of the
 representations $(V_i,KG_i)$. For every projection $\pi_i:(V,G)\to
 (V_i,G_i)$ we have epimorphism $\nu_i:(V,KG)\to
 (V_i,KG_i)$. All epimorphisms $\nu_i$ can be extended up to
 homomorphism $\mu: (V,KG)\to (V,A)$. Let $\mu(KG)=B\subset A$.
 Since $(V,A)\in \omega\Sigma$ and $(V,B)$ is a subrepresentation of $(V,A)$
 then  $(V,B)\in\omega\Sigma$. Hence  $(V,KG)\in
 \omega\Sigma$. Thus $(V,G)\in\eta\Sigma$ and the class
 $\eta\Sigma$ is a variety of representations of groups.
} 

Now we describe identities of the variety $\eta\Sigma$. Let
$\mathbb F$ be a free associative algebra over a field $K$ and let
$F$ be a free group.

\begin{defin}\label{Def_PolynIdentity}
An element $f\in\mathbb F$ is called a polynomial identity of a
representation $(V,G)$ if for any homomorphism of algebras $\mu :
\mathbb F\to KF$ the element $u=\mu(f)\in KF$ is an identity of the
representation $(V,G)$.
\end{defin}

{\bf Remark.} It is obvious, that an element $f\in\mathbb F$ is a
polynomial identity of a representation $(V,G)$ if and only if for
any homomorphism $\mu : \mathbb F\to KG$ the element $\mu(f)$
annihilates the module $V$.

Denote by $I_{\Sigma}(KF)$ the verbal ideal of $KF$ for the
variety $\Sigma$. We have the following description of the ideal
of identities of the variety $\eta\Sigma$.

\begin{prop}[\cite{PV}]\label{Prop_4_2}
The ideal $I_{\Sigma}(KF)$ is the ideal of identities of the
variety $\eta\Sigma$.
\end{prop}

\dok{
 Let $v_i$ be identities of the variety of associative algebras
 $\Sigma$, $i\in I$. The condition $(V,G)\in \eta\Sigma$
 is
 equivalent to the following one: for any
 homomorphism $\mu : \mathbb F\to KG$ the elements $\mu(v_i)$ annihilate
 the module $V$. Using the remark above we have that the
 representation $(V,G)$ belongs to $\eta\Sigma$ if and only if all
 $v_i$ are polynomial identities of $(V,G)$.

 Thus $(V,G)\in \eta\Sigma$ if and only if the images of all $v_i$
 in the algebra $KF$ are  polynomials identities of the
 representation $(V,G)$. These images generate the ideal
 $I_{\Sigma}(KF)$. So the class $\eta\Sigma$ is determined by the
 identities from the ideal $I_{\Sigma}(KF)$. Since this ideal is invariant
 under all endomorphisms of the group $F$ then it contains all
 identities of the class $\eta\Sigma$.
}

{\bf Remark.} The verbal ideal $I_{\Sigma}(KF)$ of the variety
$\Sigma$ is always invariant under all endomorphisms of the group
$F$, but in this case it is also closed under all endomorphisms of
the group algebra $KF$. However it can be shown  
that if an ideal in $KF$ 
is closed under all
endomorphisms of the group $F$ this
does not mean that this ideal is a homomorphic image of an ideal of identities for some
variety of algebras.

Remind that $\mathfrak M=\mathfrak M(K)$ is the semigroup of
varieties of representations of groups and $\mathfrak
M^{a}=\mathfrak M^{a}(K)$ is the semigroup of all non-singular
varieties of representations of associative algebras over a field
$K$.

Let us consider a map $\eta ': \mathfrak M^{a}\to \mathfrak M$ defined by the
rule $\eta '\mathfrak Y =\eta\omega^{-1}\mathfrak Y $ for every
$\mathfrak Y\in \mathfrak M^{a}$. We have the following

\begin{prop}[\cite{PV}]\label{Prop_4_3} 
The map $\eta '$ is a homomorphism of the semigroup $\mathfrak
M^{a}$ to the semigroup $\mathfrak M$.
\end{prop}

\dok{
 Let $\mathfrak Y_1$ and $\mathfrak Y_2$ be elements of the semigroup $\mathfrak
 M^{a}$ and let $\mathfrak Y=\mathfrak Y_1\mathfrak Y_2$. Let $\Sigma_1=\omega^{-1}\mathfrak
 Y_1$, $\Sigma_2=\omega^{-1}\mathfrak Y_2$ and $\Sigma=\omega^{-1}\mathfrak
 Y$. Then we have $\eta '\mathfrak Y_1=\eta\Sigma_1$, $\eta '\mathfrak
 Y_2=\eta\Sigma_2$ and $\eta '\mathfrak Y=\eta\Sigma$.

 Let $I_1$ and $I_2$ be $T$-ideals in $\mathbb F$ of the varieties
 $\Sigma_1$ and $\Sigma_2$ accordingly. Then $I=I_2I_1$ is the
 $T$-ideal of the variety $\Sigma$. Since the algebra $\mathbb F$ and the group $F$ have
 a countable rank then there is an epimorphism $\mu :\mathbb F \to
 KF$ and  $\mu(I_1)=I_{\Sigma_{1}}(KF)$, $\mu(I_2)=I_{\Sigma_{2}}(KF)$,
 $\mu(I)=I_{\Sigma}(KF)$, moreover
 $$
  \mu(I)=\mu(I_2I_1)=\mu(I_2)\mu(I_1)=I_{\Sigma_{2}}(KF)\cdot I_{\Sigma_{1}}(KF).
 $$
 The verbal ideal $I_{\Sigma_{2}}(KF)\cdot I_{\Sigma_{1}}(KF)$
 determines the variety $\eta\Sigma_{1}\cdot\eta\Sigma_{2}$. Since
 $I_{\Sigma_{2}}(KF)\cdot I_{\Sigma_{1}}(KF)=I_{\Sigma}(KF)$ then
 the variety $\eta\Sigma_{1}\cdot\eta\Sigma_{2}$ is determined by
 the ideal $I_{\Sigma}(KF)$. Hence
 $\eta\Sigma_{1}\cdot\eta\Sigma_{2}=\eta\Sigma$.
} 

We can give the following

\begin{defin}\label{Def_PolynVariety}
A variety of representations of groups $\mathfrak X$ is called a
polynomial variety if there exists a variety of associative
algebras $\Sigma$ such that $\mathfrak X=\eta\Sigma$.
\end{defin}

Let us give some examples of polynomial varieties.


\bigskip
 {\bf Examples.}
\begin{enumerate}
 \item
  Let $\Sigma$ be the variety of all commutative algebras over a field $K$
  and let $\mathfrak A$ be a variety of abelian groups. Then  $\omega \mathfrak A=\eta
  \Sigma$. 

 \item
  Consider the variety $(\omega\mathfrak A)^{n}$. Using the
  homomorphism $\eta ':\mathfrak M^{a}\to \mathfrak M$ we have
  $\eta '(\omega \Sigma)=\eta \Sigma=\omega\mathfrak A$. Further
  $$
    \eta '((\omega \Sigma)^{n})=\eta(\omega^{-1}(\omega \Sigma)^{n}),
  $$
  $$
   \eta '((\omega \Sigma)^{n})=(\omega \Sigma)^{n}=(\mathfrak S\times \mathfrak
   A)^{n}=\mathfrak S^{n}\times \mathfrak A,
  $$
  (see section \ref{SubS1.3_AlgOfVar}). A variety of algebras
  $\Sigma_1=\omega^{-1}((\omega \Sigma)^{n})$
  is a variety of algebras with the nilpotent
  derived ideal of class $n$. Thus we have $\mathfrak S^{n}\times \mathfrak
  A=\eta\Sigma_1$. So the variety $\mathfrak S^{n}\times \mathfrak
  A$ is  polynomial. Remind that the variety $\mathfrak S^{n}$ is a variety
generated by the representation $(K^{n},UT_{n}(K))$, and the
variety $(\omega \mathfrak A)^{n}$ is generated by the
representation $(K^{n},T_{n}(K))$.

\end{enumerate}

Proposition~\ref{Prop_4_3} implies:

\begin{cor}[\cite{PV}]\label{Cor_4_1} 
All polynomial varieties of representations of groups constitute
a subsemigroup in the semigroup $\mathfrak M$. $\qquad \square$
\end{cor}

Remind that if $K$ is a field then the semigroup $\mathfrak M$ of
all varieties of representations of groups is free. So we have the
following

\begin{prob}
Is the semigroup of all polynomial varieties free?
\end{prob}

To solve this problem we study some properties of polynomial
varieties of representations of groups. In particular, we are
looking for interesting examples of polynomial and non-polynomial
varieties.

For any representation $(V,G)$ the variety $\mbox{Var}(V,G)$ is the
variety generated by the representation $(V,G)$. We denote by
$\widetilde{\mbox{Var}}(V,G)$ the minimal polynomial variety which
contains the representation $(V,G)$, in other words the variety
$\widetilde{\mbox{Var}}(V,G)$ is the minimal polynomial variety
containing the variety $\mbox{Var}(V,G)$. It is not difficult to
see that $\widetilde{\mbox{Var}}(V,G)=\eta\Sigma$, where
$\Sigma=\mbox{Var}(\overline{KG})$ and
$\overline{KG}=KG/\mbox{Ker}(V,KG)$.

For any variety of representations of groups $\mathfrak X$ we
denote by $\widetilde{\mathfrak X}$ the minimal polynomial variety
containing the variety $\mathfrak X$.



%

%
%

Note that if $\mathfrak Y=\mbox{Var}(V,\overline{KG})$ is a
variety of representations of algebras and
$\Sigma=\mbox{Var}(\overline{KG})$ then
$\Sigma=\omega^{-1}\mathfrak Y$, (remind that $\omega^{-1}\mathfrak
Y=\overrightarrow{\mathfrak Y}$ is the class of all algebras
admitting a faithful representation in the variety $\mathfrak Y$).
Indeed, since the representation $(V,\overline{KG})$ is faithful
then $\overline{KG}\in\omega^{-1}\mathfrak Y$ and so
$\Sigma\subset \omega^{-1}\mathfrak Y$. On the other hand, if $A$
is an algebra from $\omega^{-1}\mathfrak Y$ then there exists a
$K$-module $U$ such that $(U,A)$ is a faithful representation from
the variety $\mathfrak Y$. If $y\circ w(x_1,\dots ,x_n)\equiv 0 $
is an identity in $\mathfrak Y$ then $w(x_1,\dots ,x_n)\equiv 0$
is an identity of the algebras $\overline{KG}$ and  $A$, since
$(V,\overline{KG})$ and $(U,A)$ are faithful representations. So
$A\in \Sigma=\mbox{Var}(\overline{KG})$. Thus
$\Sigma=\omega^{-1}\mathfrak Y$.

For polynomial varieties we have the property~\ref{Prop_Example}
which show that the proposition~\ref{Prop_1-1} is not true for
polynomial varieties. 

\begin{prop}\label{Prop_Example}
There are representations  $(V_1,G_1)$ and $(V_2,G_2)$ such that
$$
\widetilde{\mbox{Var}}((V_1,G_1)\bigtriangledown (V_2,G_2))\ne
\widetilde{\mbox{Var}}(V_1,G_1)\cdot
\widetilde{\mbox{Var}}(V_2,G_2).
$$
\end{prop}

\dok{
 Let $K$ be a field and let $(V_1,G_1)=(V_2,G_2)=(K,1)$.
 Then
 $$
  \widetilde{\mbox{Var}}((K,1)\bigtriangledown
  (K,1))=\widetilde{\mbox{Var}}(K^{2},\mbox{UT}_2(K))=\eta\Sigma,
 $$
 where $\Sigma=\mbox{Var}A$ and
 $A$ is the algebra of matrices of the form
 $
 \left(
 \begin{array}{cc}
 \alpha & \mbox{End} K \\
 0      & \alpha \\
 \end{array}
 \right)$,
 $\alpha\in K$. Note that the variety $\Sigma$ is the  variety
 of commutative algebras.

 In the section~\ref{S4_VarOfReprGrVarAlg} we considered the homomorphism
 $\eta': \mathfrak M^{a}\to \mathfrak M$ of the semigroup $\mathfrak M^{a}$ of
 varieties of representations of algebras to the semigroup $\mathfrak M$
 of varieties of representations of groups. Remind that if $\mathfrak Y\in \mathfrak
 M^{a}$ then $\eta'\mathfrak Y=\eta\omega^{-1}\mathfrak Y$.

 We have the following equalities:
 $$
  \eta'\mbox{Var}(V,\overline{KG})=\eta\omega^{-1}\mbox{Var}(V,\overline{KG})=
  \eta\mbox{Var}(\overline{KG})=\widetilde{\mbox{Var}}(V,G).
 $$

 Along with the triangular product of representations of groups we
 can consider the triangular product of representations of
 algebras.  It is known \cite{Kaljulaid} that for varieties of
 representations of associative algebras over a field there is the following
 property:
  $$
   \mbox{Var}((V_1,A_1)\bigtriangledown
   (V_2,A_2))=\mbox{Var}(V_1,A_1)\cdot \mbox{Var}(V_2,A_2).
  $$

 We denote by $K\cdot 1$ the group algebra over $K$ of the trivial
 group.
 So we have
 $$
 \widetilde{\mbox{Var}}((K,1))\cdot \widetilde{\mbox{Var}}((K,1))=
 \eta'\mbox{Var}(K,K\cdot 1)\cdot
  \eta'\mbox{Var}(K,K\cdot 1)=
 $$
 $$
 =\eta'\Big(\mbox{Var}(K,K\cdot 1)\cdot \mbox{Var}(K,K\cdot 1)
  \Big)=
  \eta'\Big( \mbox{Var}((K,K\cdot 1)\bigtriangledown (K,K\cdot
  1))\Big)=
 $$
 $$
 =\eta'\Big( \mbox{Var}(K^2,A_1)
 \Big),
 $$
 where $A_1$ is an algebra of matrices of the form
 $\left(
 \begin{array}{cc}
  K\cdot 1 & \mbox{End} K\\
  0        & K\cdot 1 \\
 \end{array}
 \right)$ and this algebra is not commutative.

 It is known \cite{PV} 
  that a unitriangular product of representations of groups is
  faithful if and only if all  factors are faithful.
 Since the representation $(K,K\cdot 1)$ is faithful then the
 representation $(K^2,A_1)$ is also faithful and we have
 $$
 \eta'\Big( \mbox{Var}(K^2,A_1)
 \Big)=\eta\mbox{Var}(A_1),
 $$
 where $\mbox{Var}(A_1)$ is not commutative variety. So we have
 that
 $$
 \widetilde{\mbox{Var}}((K,1)\bigtriangledown (K,1))\ne
\widetilde{\mbox{Var}}(K,1)\cdot \widetilde{\mbox{Var}}(K,1).
 $$
}

We have following property for the variety $\mathfrak S^{n}$.

\begin{prop}\label{Prop_4_4}
 The variety of representations of groups $\mathfrak S^{n}$ is not a
 polynomial variety.
\end{prop}

\dok{
 To prove the proposition it is enough to show that the equality $\mathfrak
 S^{n}=\eta\Sigma$ is impossible for any variety $\Sigma$.

 Assume that there is a variety of associative algebras $\Sigma$ such that $\mathfrak
 S^{n}=\eta\Sigma$. The variety $\mathfrak S^{n}$ contains the representation $(K^{n},G)$,
 where $G=UT_{n}(K)$ is the group of
 unitriangular matrices and the group $G$ acts in a natural way on $K^{n}$.
 The representation $(K^{n},G)$ belongs to the variety $\eta\Sigma$ if
 the linear span $\langle G \rangle_K$ of group $G$ in the algebra
 $\mbox{End}(K^{n})=M_n(K)$ is an algebra from the variety  $\Sigma$.
 Note that the representation $(K^{n},G)$ is faithful, so $\overline
 G=G$.


 The algebra $\langle G \rangle_K$ is an algebra of matrices of the form:
 $$
 \langle G \rangle_K=
 \left\{
 \left(
 \begin{array}{ccc}
 \alpha & * & *\\
  0     &\ddots & *\\
   0    & \dots & \alpha\\
 \end{array}
 \right),
 \alpha \in K
 \right\}
 $$
 Let $G_1$ be a group of all invertible elements of the algebra $\langle G \rangle_K$.
 The group $G_1$ acts faithfully on the module $K^{n}$ and the
 linear span $\langle G_1 \rangle_K$ of  $G_1$ in the algebra
 $\mbox{End}(K^{n})$ coincides with $\langle G \rangle_K$ which is an algebra
 from the variety $\Sigma$.
 %
 %
 So the representation $(K^{n},G_1)$ belongs to the variety $\eta\Sigma=\mathfrak S^{n}$.
 But it is not true. Thus $\mathfrak S^{n}\ne \eta\Sigma$.
} 


We have proved that $\mathfrak S^{n}$ is not a polynomial variety.
From the proof of the proposition follows that the minimal
polynomial variety which contains $\mathfrak S^{n}$ is the variety
$\widetilde{\mathfrak S^{n}}=\eta\Sigma$, where $\Sigma$ is the
variety of algebras generated by the algebra of matrices of the
following from:
 $$
  \left\{
 \left(
 \begin{array}{ccc}
 \alpha & * & *\\
  0     &\ddots & *\\
   0    & \dots & \alpha\\
 \end{array}
 \right),
 \alpha \in K
 \right\}
 $$

Remind that the variety $\omega \mathfrak A$ is a polynomial
variety (see examples above) and $\omega \mathfrak A=\eta\Sigma$,
where $\Sigma$ is a variety of commutative algebras. Moreover, the
minimal polynomial variety $\widetilde{\mathfrak S}$ which
contains $\mathfrak S$ coincides with $\omega \mathfrak A$, that
is $\widetilde{\mathfrak S}=\omega \mathfrak A$, since the variety
$\mathfrak S=\mbox{Var}(K,1)$ and the algebra $\langle 1
\rangle_K$ generates the variety of commutative algebras.

 Let $\mathfrak X$ be a variety of representations of groups
over a field $K$. Let $I_{\mathfrak X}(KF)\triangleleft KF$ be the
ideal of identities of the variety $\mathfrak X$ and $I_{\mathfrak
X}(KG)$ be the verbal ideal of the group algebra $KG$ . Let
$\Sigma $ be a variety of associative algebras over a field $K$.

\begin{defin}
Denote by $\Sigma * \mathfrak X$ a class of representations of
groups over a field $K$ such that a representation $(V,G)$ belongs
to $\Sigma * \mathfrak X$ if the representation $(V,I_{\mathfrak
X}(KG))$ lies in $\omega \Sigma$.
\end{defin}

In other words, the representation $(V,G)$ belongs to $\Sigma
* \mathfrak X$ if and only if there is a right ideal $J$ in the
algebra $KG$ such that $(V,J)\in\omega\Sigma$ and $(KG/J,
G)\in\mathfrak X$.

We have

\begin{prop}[\cite{PV}]\label{Prop_4_5} 
For any variety of associative algebras $\Sigma$ and any variety
of representations of groups $\mathfrak X$ the class $\Sigma
* \mathfrak X$ forms a variety of representations of groups. $\qquad
\square$
\end{prop}

Let ${\bf N}_{n}$ be a variety of nilpotent algebras of nilpotency
class $n$ (that is ${\bf N}_{n}$ is the variety satisfying the
identity $x_1x_2\dots x_n=0$). Let $\mathfrak X$ be a variety of
representations of groups and let $I_{\mathfrak X}(KF)$ be the
ideal of identities for the variety $\mathfrak X$.
We have the following

\begin{prop}
${\bf N}_n * \mathfrak X=\mathfrak X^{n}$.
\end{prop}

\dok{ A representation $(V,G)$ belongs to the variety ${\bf N}_{n}
* \mathfrak X$ if and only if the algebra $\overline{ I_{\mathfrak
X}(KG)}=I_{\mathfrak X}(KG)/\mbox{Ker}(V,I_{\mathfrak X}(KG))$ is
an algebra from the variety ${\bf N}_{n}$, that is
$\big(\overline{I_{\mathfrak X}(KG)}\big)^{n}=0$. So for every
$v\in V$ we have $v\circ \big(\overline{I_{\mathfrak
X}(KG)}\big)^{n}=0$. Since $\big(I_{\mathfrak
X}(KG)\big)^{n}=I_{\mathfrak X^{n}}(KG)$ then
$\big(\overline{I_{\mathfrak
X}(KG)}\big)^{n}=\overline{I_{\mathfrak X^{n}}(KG)}$ and $v\circ
\big(\overline{I_{\mathfrak X}(KG)}\big)^{n}=0$ if and only if
$v\circ I_{\mathfrak X^{n}}(KG)=0$. So the representation
$(V,G)\in \mathfrak X^{n}$.

Thus ${\bf N}_n * \mathfrak X=\mathfrak X^{n}$.

} 

%

Define a PI-variety of representations of groups 

\begin{defin}\label{Def_PI-Variety}
A variety of representations of groups $\mathfrak X$ is called a
PI-variety if every representation from $\mathfrak X$ satisfies
some non-trivial polynomial identity.
\end{defin}

There is the following

\begin{prop}[\cite{PV}]\label{Prop_4_6} 
If a variety $\mathfrak X$ is generated by a representation
$(V,G)$ satisfying  a polynomial identity $f$ then all
representations from $\mathfrak X$ satisfy this identity.
\end{prop}

\dok{
 Let $\mathfrak X$ be a variety generated by the representation
$(V,G)$ and let $(V,G)$ satisfies  a polynomial identity $f$.
 Let $\Sigma$ be a variety of algebras determined by the identity
 $f$. Then $(V,G)$ belongs to a variety $\eta\Sigma$. Hence $\mbox{Var}(V,G)=\mathfrak
 X$ is contained in the variety $\eta\Sigma$. So all representations from $\mathfrak X$
 satisfy the polynomial identity $f$.
} 

\begin{cor}[\cite{PV}]\label{Cor_4_2} 
If $\mathfrak X$ is a polynomial variety then there is a
non-trivial polynomial identity such that all representations from
$\mathfrak X$ satisfy this identity. $\qquad \square$
\end{cor}

{\bf Remark.} Let $\mathbb F$ be a free associative algebra over a
field $K$ with free generates $z_1,z_2,\dots $ and let $F$ be a
free group with free generates $x_1,x_2,\dots $. Now let
$f=f(z_1,\dots ,z_n)\in \mathbb F$ be a multilinear polynomial.
Consider a map $\mathbb F\to F$ by the rule $z_i\to x_i$. The
polynomial $f=f(z_1,\dots z_n)$ is an identity of a representation
$(V,G)$ if and only if the representation $(V,G)$ satisfies the
identity $f'=f(x_1,\dots ,x_n)$.

It is well known that every $PI$-algebra satisfies a multilinear
identity. Hence every $PI$-variety of representations of groups
satisfies  some non-trivial multilinear identity.  The remark
above allows us to consider this identity as an element of the group
algebra $KF$.

If $K$ is a field of zero characteristic then every polynomial
identity are equivalent to some multilinear identity. Thus every
polynomial varieties of representations of groups over a field of
zero characteristic can be determined by the  multilinear
identities.

\section{General view on dimension subgroups}\label{S2_DimSubgr}

Recall that the classical dimension subgroup $D_{n}(\mathbb Z,G)$
is a subgroup of the group $G$ over $\mathbb Z$ is defined as
follows: $D_{n}(\mathbb Z,G)$ is the set of all elements $g\in G$
such that $(g-1)$ belongs to the ideal $(\Delta _{G})^{n}\subset
\mathbb Z G$, where $\Delta_{G}$ is the augmentation ideal of
$\mathbb Z G$. That is
$$
D_{n}(K,G)=(1+(\Delta _{G})^{n})\cap G.
$$

The exposition of the classical dimension subgroup, various
problems concerning dimension subgroups and book review  are
contained in the book of R.~Mikhailov and
I.B.S.~Passi~\cite{MikhPassi_book}. In particular, we would note
the papers by W.~Magnus~\cite{Magnus}, O.~Gr\" un~\cite{Grun},
N.~Gupta~\cite{Gupta}, I.B.S.~Passi~\cite{Passi},
E.~Rips~\cite{Rips}, E.~Witt~\cite{Witt}. We are interested to
consider the dimension subgroup from more general point of view.
Remind that the n-th power  $(\Delta _{G})^{n}\subset KG$ of the
augmentation ideal $(\Delta _{G})$ is  the ideal of identities of
the variety of representations of groups $\mathfrak S^{n}$. So
this insight yields a new approach to dimension subgroups.

\subsection{Dimension subgroup for a variety of representations of groups}

Let $\mathfrak X$ be a variety of representations of groups. We
consider  representations $(V,G)$ of the group $G$
contained in the variety $\mathfrak X$. We define a subgroup
$D_{\mathfrak X}(G)$ of the group $G$ as follows:

\begin{defin}\label{dim_subg_var}
A dimension subgroup $D_{\mathfrak X}(G)$ of a group $G$ for a
variety of representations of groups $\mathfrak X$ is the
intersection of kernels of all representations $(V,G)\in \mathfrak X$. 
\end{defin}

Denote by $\overrightarrow{ {\mathfrak X}}$ the class of groups
admitting faithful representation in the variety $\mathfrak X$.
This class is a quasi-variety of groups and the corresponding
verbal subgroup coincides with the dimension subgroup. In
particular, if $K$ is a field of zero characteristic then
$\overrightarrow{\mathfrak S^{n}}$ is the quasi-variety of all
torsion free nilpotent groups of the class $(n-1)$ \cite{Malcev2}.
For the quasi-variety $\overrightarrow{ {\mathfrak S^{n}}}$ there
is also the problem of B.~Plotkin~\cite{PV} which asks: has the
quasi-variety $\overrightarrow{ {\mathfrak S^{n}}}$ a finite basis
of quasi-identities? For the quasi-variety $\overrightarrow{
{\mathfrak S^{4}}}$ this problem is solved  by R.~Mikhailov and
I.B.S.~Passi~\cite{MikhPassi_D4} (see also \cite{MikhPassi_book}).
They proved that this quasi-variety is not finitely based.

Let $I_{\mathfrak X}(KF)$ be the ideal of identities of the group
algebra $KF$ for a variety $\mathfrak X$ (here $F$ is a free
group) and let $I_{\mathfrak X}(KG)$ be the verbal ideal of $KG$
for the variety $\mathfrak X$. The representation
$(KG/I_{\mathfrak X}(KG),G)$ is the free representation of the
group $G$ in the variety $\mathfrak X$. The kernel of this
representation is a normal subgroup $H=(1+I_{\mathfrak X}(KG))\cap
G$ consisting of all $g\in G$ such that the element $(g-1)$
belongs to the ideal $I_{\mathfrak X}(KG)$.

We have the following proposition.

\begin{prop}[\cite{PV}]\label{Prop_KernelOfUnivLin} 
The kernel of the  free representation $(KG/I_{\mathfrak
X}(KG),G)$ coincides with the dimension subgroup $D_{\mathfrak
X}(G)$ of the group $G$ for the variety $\mathfrak X$.
\end{prop}

\dok{
 By the definition of $D_{\mathfrak X}(G)$ we have  $D_{\mathfrak
X}(G)\subseteq H$, where $H=(1+I_{\mathfrak X}(KG))\cap G$ is the
kernel of the  free representation $(KG/I_{\mathfrak X}(KG),G)$.
We will show that $H \subseteq D_{\mathfrak X}(G)$.

 Let $(V,G)$ be a representation of the group $G$ such that $(V,G)\in \mathfrak X$, 
 let $v$ be an element in $V$. The map $g\to v\circ g$ induces the
 homomorphism of modules: $KG\to V$,  and this homomorphism is
 permutable with the action of the group $G$. Since the ideal $I_{\mathfrak
 X}(KG)$ annihilates the module $V$ then there is a homomorphism $\nu: KG/I_{\mathfrak X}(KG)\to
 V$. Let $g$ be an arbitrary element from the kernel $H$ and let $\overline 1=1+I_{\mathfrak X}(G)$. 
 Then for all $v\in V$
 we have:
 $$
  v\circ g=(\overline 1)^{\nu}\circ g=(\overline 1\circ
  g)^{\nu}=(\overline 1)^{\nu}=v.
 $$
 So any element $g\in H$ belongs to the kernel of the representation
 $(V,G)$. It means that the kernel of the representation $(KG/I_{\mathfrak X}(KG), G)$ is
 contained in the kernel of every representation $(V,G)$ form the variety $\mathfrak X$.
 Thus $D_{\mathfrak X}(G)$ coincides with the kernel of the  free representation $(KG/I_{\mathfrak
X}(KG),G)$.
} 

Since the ideal $(\Delta _{G})^{n}$ is the ideal of identities of
the variety $\mathfrak S^{n}$
the classical dimension subgroup $D_{n}(K,G)$ is the dimension
subgroup of the group $G$  for the variety $\mathfrak X=\mathfrak
S^{n}$, that is $D_{n}(K,G)=D_{\mathfrak S^{n}}(G)$. So we have
that the classical dimension subgroup connects with the dimension
subgroup for the variety of representations of groups $\mathfrak
S^{n}$.

It is interesting to consider the following

\begin{prob}
Describe the dimension subgroups for the polynomial varieties
$\widetilde{\mathfrak S^{n}}$ and $(\widetilde{\mathfrak
S})^{n}=(\omega \mathfrak A)^{n}$, respectively.
\end{prob}

We have a partial solution of this problem for the free group $F$
and the variety $\widetilde{\mathfrak S^{n}}$, see
Theorem~\ref{Th_EqualDimSubgrFree} below.

 Remind that the variety $\mathfrak S^{n}$ is the variety
generated by the representation $(K^{n},UT_{n}(K))$, and the
variety $(\widetilde{\mathfrak S})^{n}=(\omega \mathfrak A)^{n}$
is generated by the representation $(K^{n},T_{n}(K))$.

In the section~\ref{SubS5.3}
we consider the concept of a dimension subgroup for variety of
associative algebras and some relations between these subgroups.

\subsection{Linearization of groups and
$PI$-groups}\label{S5_Linearization}

\paragraph{5.2.1}

Let $A$ be an associative algebra with unit over a field $K$ and
let $G$ be a group.

\begin{defin}\label{def_Lineariz}
An arbitrary  homomorphism $\rho$: $G\to A\ $ of a group $G$ into
the group of invertible elements of an algebra $A$ is called a
linearization of the group $G$ in the algebra $A$ over a field
$K$.
\end{defin}

For a given group $G$ we  consider the category of all
$K$-linearizations of  $G$ over a field $K$. The objects
of this category are linearizations $\rho_i :G\to A_i$ of the
group $G$ in the $K$-algebras $A_i$. Let $\rho_1$ and $\rho_2$ be
two linearizations of the group $G$ in the algebras $A_1$ and
$A_2$, accordingly.
%
%
Then there is a homomorphism $\mu$ of the algebra $A_1$ to the
algebra $A_2$  such that 

\begin{picture}(100,60)\
\put(55, 20){\vector(2,1){40}} \put(40,20){$G$} \put(100,
35){$A_1$} \put(70,35){$\rho_1 $}%
\put(55, 20){\vector(2,-1){40}} \put(100,-5){$A_2$} \put(70,-5){$\rho_2$}%
\put(105, 30){\vector(0,-1){20}} \put(110, 15){$\mu$}
\end{picture}

\bigskip

%
%

 So the morphisms of this category are commutative diagrams  above.

This category possesses a universal initial object. It is
identical embedding $G\stackrel{id}{\longrightarrow } KG$  of the
group $G$ into the group algebra $KG$. For every linearization
$\rho: G\to A$ we have a homomorphism $\mu :KG\to A$ such that the
following diagram is commutative

\begin{picture}(100,60)\
\put(55, 35){\vector(1,0){40}} \put(40,30){$G$}\put(70,40){$_{id}$} \put(95,30){$KG$}%
\put(55, 30){\vector(2,-1){40}} \put(95,0){$A$} \put(70,10){$\rho$}%
\put(100, 25){\vector(0,-1){15}} \put(105, 15){$\mu$}
\end{picture}


\begin{defin}\label{Def_FaithfulLinear}
A linearization $\rho :G\to A$ is called  faithful 
if the map $\rho$ is an injective homomorphism.
\end{defin}

The linearization $G\to KG$ is faithful. In  general case, if
$\rho : G\to A$ is a linearization of a group $G$ and
$G^{\rho}=G/\mbox{Ker}\rho$ then the linearization $\overline\rho
:G^{\rho}\to A$ is faithful.

If a linearization $\rho :G\to A$ is a faithful linearization then
the group $G$ can be considered as a subgroup of the group of
invertible elements of the algebra $A$.

Let $\Sigma $ be a non-trivial variety of associative algebras
over a fixed field $K$. We will consider a linearization $\rho
:G\to A$, where $A$ is a $K$-algebra from the variety $\Sigma$.

We also have the category of such linearizations and its initial
object is a linearization $\rho :G\to KG/I_{\Sigma}(KG)$, where
$I_{\Sigma}(KG)$ is the verbal ideal in $KG$ for the variety
$\Sigma$. This linearization can be not faithful.

\begin{defin}\label{def_PI-gr}
A group $G$ is called a $PI$-group if there is a faithful
linearization $\rho :G\to A$, where $A$ is an algebra from some
non-trivial variety of associative algebras $\Sigma$.
\end{defin}

In more detail we consider $PI$-group in \cite{AladovaPlotkin} where
we solve for these groups some Burnside-type problems.

Note that,  if a group $G$ admits a faithful linearization in an algebra
$A$ from a variety $\Sigma$ then the linearization $G\to
KG/I_{\Sigma}(KG)$ is also faithful.

\paragraph{5.2.2}

Let us return to  dimension subgroups.

Let $\rho: G\to A$ be a linearization of a group $G$ in an algebra
$A$ from a variety $\Sigma$. In particular we can consider
universal linearization $G\to KG/I_{\Sigma}(KG)$, where
$I_{\Sigma}(KG)$ is the  verbal ideal in $KG$ of the variety
$\Sigma$. It is clear that the kernel of every homomorphism
$\rho:G\to A$ contains the kernel of universal linarization. 
Observe that the kernel of the universal linarization is the set
of all $g\in G$ such that the element $(g-1)\in KG$ belongs to the
ideal $I_{\Sigma}(KG)$. This kernel is denoted by $D_{\Sigma}(G)$. So
we have
$$
 D_{\Sigma}(G)=\{ \ g\in G \mid\ g-1\in I_{\Sigma}(KG) \ \}.
$$

\begin{defin}\label{Def_DimSubgr}
The kernel $D_{\Sigma}(G)$ of universal linearization $G\to
KG/I_{\Sigma}(KG)$ is called a dimension subgroup of the group $G$
for the variety of algebras $\Sigma$.
\end{defin}

Note that a group $G$ is a $PI$-group with given variety $\Sigma$
if and only if $D_{\Sigma}(G)=1$.

Let $\mathfrak D(\Sigma)$ be a class of all $PI$-groups with the given
$\Sigma$. Here we have a theorem similar to  the theorem  for the
varieties of representations of groups. The latter one states that a class of
all groups $G$ admitting a faithful representation $(V,G)$ from
some variety of representations of groups $\mathfrak X$ forms a
quasivariety of groups~\cite{PV}.


\begin{teo}\label{Th_Quasivar}
The class $\mathfrak D(\Sigma)$  is a quasivariety of groups.
\end{teo}

\dok{ To prove the theorem it is sufficient to show \cite{Malcev}
that the class $\mathfrak D(\Sigma)$ contains the trivial group
and  is closed
under subgroups and filtered products. It is obvious that for 
class
$\mathfrak D(\Sigma)$ first two conditions hold true. We need
to prove that the class $\mathfrak D(\Sigma)$ is closed under the
filtered products.

Let $I$ be a non-empty set, let $D$ be a filter on the set $I$.
Let $G_{\alpha}$ be groups from the class $\mathfrak D(\Sigma)$,
let $\varphi_{\alpha}$ : $G_{\alpha} \to A_{\alpha}$ be faithful
linearizations of groups $G_{\alpha}$ in the algebras $A_{\alpha}$
from the variety $\Sigma$.

Let $A=\prod_{\alpha} A_{\alpha}$ be the Cartesian product of the
algebras $A_{\alpha}$. Let $A/J$ be the filtered product of these
algebras, where $J$ is the ideal of $A$ such that
$$
 a\in J
\Leftrightarrow \{\ \alpha\in I \ \mid\ a(\alpha)=0 \}\subset D.
$$


Since all $A_{\alpha}$ are PI-algebras from the variety $\Sigma$
then the algebra $A$ and the algebra $A/J$ belongs to the variety
$\Sigma$.

Let $G=\prod_{\alpha} G_{\alpha}$ be the Cartesian product of the
groups $G_{\alpha}$.
 Let $G/H$ be the filtered product of the groups
$G_{\alpha}$,  where $H$ is the subgroup of $G$ such that
$$
 h\in H\Leftrightarrow \{\ \alpha\in I \ \mid\ h(\alpha)=1 \}\subset D.
$$



Consider the linearization  $\psi\ :G\to A/J$ defined by the rule $f\to
f'+J$ where $f'(\alpha)=\varphi_{\alpha}(f(\alpha))$ for all
$\alpha \in J$, $f\in G$. It is easy to see that $\psi$ is a
homomorphism.


Consider the kernel of the homomorphism $\psi$:
$$
\mbox{Ker}(\psi)=\{\ g\ \mid \ \psi(g)=1+J\ \}.
$$

For any element $g\in \mbox{Ker}(\psi)$ we have $ \psi(g)(\alpha
)=\varphi_{\alpha}(g(\alpha))+J(\alpha)=1(\alpha)+J(\alpha)=1(\alpha)$,
so $\varphi_{\alpha}(g(\alpha))=1(\alpha)$ and $g\in H$. Thus
$\mbox{Ker}(\psi)=H$.


Then $G/H\cong \mbox{Im }\psi \subset A/I$ and the group $G/H$
admits the faithful linearization in the algebras from the variety
$\Sigma$. Thus the class $\mathfrak D(\Sigma)$ is closed under the filtered
products.

So the class $\mathfrak D(\Sigma)$ is a quasivariety.
} 

\subsection{Dimension subgroup for varieties of representations of
groups and for varieties of associative algebras}\label{SubS5.3}

\paragraph{5.3.1}

Now we will compare two approaches (see Definition \ref{dim_subg_var}   and Definition  \ref{Def_DimSubgr}) to the notion of a dimension subgroup.


 Let $\mathfrak X$ be a variety of representations of groups
and let $I_{\mathfrak X}(KG)$ be the verbal ideal in $KG$ for the
variety $\mathfrak X$. Let $\Sigma$ be a variety of associative
algebras over a field $K$ and let $I_{\Sigma}(KG)$ be the verbal
ideal in $KG$ for the variety $\Sigma$. In the section {\bf
\ref{S4_VarOfReprGrVarAlg}} (see Proposition {\bf \ref{Prop_4_2}})
it was proved that if $\mathfrak X=\eta\Sigma$ then $I_{\mathfrak
X}(KF)=I_{\Sigma}(KF)$, where $F$ is the free group. 

Hence, here we have the following

\begin{teo}\label{Th_EqualDimSubgr}
If $\mathfrak X=\eta\Sigma$ is a polynomial variety  then the
dimension subgroup $D_{\mathfrak X}(G)$ coincides with the
dimension subgroup $D_{\Sigma}(G)$. $\qquad \square$
\end{teo}

{\bf Remark.} Remind that there is an equality $\mathfrak S^{n}
\times \mathfrak A=(\omega \mathfrak A)^{n}=\eta\Sigma$ (see
section {\bf\ref{S4_VarOfReprGrVarAlg}}), where $\mathfrak A$ is
the variety of abelian groups and  $\Sigma$ is the variety of
associative algebras satisfying the identity
$$
(x_1y_1-y_1x_1)\dots (x_ny_n-y_nx_n)\equiv 0.
$$
So, we have the following equality
 $$
 D_{\mathfrak S^{n} \times \mathfrak A}=D_{\Sigma}.
 $$

Remind that $D_n(K,G)=(1+(\Delta _{G})^{n})\cap G$ is the
dimension subgroup of a group $G$ over a ring $K$, where
$\Delta_{G}$ is the augmentation ideal of $KG$, and for every
group $G$ the dimension subgroup $D_n(K,G)$ coincides with the
dimension subgroup $D_{\mathfrak S^{n}}(G)$ for the variety of
representations of groups $\mathfrak S^{n}$. For a free group $F$
we have the following

\begin{teo}\label{Th_EqualDimSubgrFree}
For a free group $F$  its dimension subgroup $D_n(K,F)$ coincides
with the dimension subgroup $D_{\widetilde{\mathfrak S^{n}}}(F)$
of the group $F$ for the polynomial variety of representations of
groups $\widetilde{\mathfrak S^{n}}$.
\end{teo}

\dok{

In the section~{\bf \ref{S4_VarOfReprGrVarAlg}} it was shown that
$\widetilde{\mathfrak S^{n}}=\eta\Sigma$, where $\Sigma$ is a
variety of algebras generated by the algebra of matrices of the
following from:
 $$
  \left\{
 \left(
 \begin{array}{ccc}
 \alpha & * & *\\
  0     &\ddots & *\\
   0    & \dots & \alpha\\
 \end{array}
 \right),
 \alpha \in K
 \right\},
 $$
and this algebra is the linear span $\langle G \rangle_{K}$ of the
group $G=UT_{n}(K)$  of unitriangular matrices  in the algebra
$\mbox{End}(K^{n})$.

 Let $H=(\langle G \rangle_{K})^{*}$ be the group of invertible
elements of the algebra $\langle G \rangle_{K}$, that is
$$
H=
 \left\{
  \left(
    \begin{array}{ccc}
     \alpha & * & *\\
       0    &\ddots & *\\
       0    & \dots & \alpha\\
    \end{array}
      \right)
      \mid \alpha\in K, \alpha\ne 0
\right\}.
$$
The group $H=\mbox{UT}_n(K)\cdot Z(T_{n}(K))=\mbox{UT}_n(K)\times
Z(\mbox{T}_{n}(K))$, where $Z(\mbox{T}_{n}(K))$ is the center of
the group $\mbox{T}_n(K)$. The group $\mbox{UT}_n(K)$ is nilpotent
of the nilpotency class $(n-1)$ and the group $Z(\mbox{T}_{n}(K))$
is abelian hence the group $H$ is nilpotent of the class $(n-1)$.

Let $F$ be a free group. Consider all possible homomorphisms of
the group $F$ to the group $H$. These homomorphisms induce
homomorphisms of the group algebra $KF$ to the group algebra
$KH=\langle G \rangle_{K}$:
$$
\psi_{i}: KF\to \langle G \rangle_{K}.
$$
Let $U_{i}$ be the kernel of the homomorphism $\psi_i$ and let
$U=\bigcap_{i}U_{i}$ be the intersection of kernels of all
possible homomorphisms $\psi_{i}$. So $U$ is the set of all
identities of the algebra $\langle G \rangle_{K}$ in the algebra
$KF$, in other words, $U$ is the verbal ideal of the variety
$\Sigma$ in the group algebra $KF$:
$$
U=I_{\Sigma}(KF).
$$
Using the proposition~\ref{Prop_4_2} we have:
$$
I_{\Sigma}(KF)=I_{\eta\Sigma}(KF)=I_{\widetilde{\mathfrak
S^{n}}}(KF).
$$

According to the Remak's theorem the algebra $KF/U$ is a
subalgebra of the Cartesian product of  $KF/U_{i}$. All
algebras $KF/U_{i}$ are subalgebras of $\langle G \rangle_{K}$. So
the algebra $KF/U$ can be embedded into a Cartesian
power of the 
algebra $\langle G
\rangle_{K}$. Denote this cartesian power by $\prod \langle G \rangle_{K}$.

Then,  $\prod H$ is the group of  all invertible elements of the
algebra $\prod \langle G \rangle_{K}$. So the group $(KF/U)^{*}$
of all invertible elements of the algebra $KF/U$ is a subgroup of
$\prod H$.

The group $H$ is nilpotent of the class $(n-1)$. Thus,  the group
$\prod H$ is also nilpotent of the class $(n-1)$. Hence the group
$(KF/U)^{*}$ is nilpotent of nilpotency class at most $(n-1)$.

Let $D_{\widetilde{\mathfrak S^{n}}}(F)$ be the dimension subgroup
of the free group $F$ for the variety $\widetilde{\mathfrak
S^{n}}$, that is
$$
D_{\widetilde{\mathfrak S^{n}}}(F)=\{ g\in F\mid (g-1)\in
U=I_{\widetilde{\mathfrak S^{n}}}(KF) \}.
$$
Since $I_{\widetilde{\mathfrak S^{n}}}(KF)=I_{\Sigma}(KF)$
then $D_{\widetilde{\mathfrak S^{n}}}(F)=D_{\Sigma}(F)$ and the
subgroup $D_{\Sigma}(F)$ is a kernel of the homomorphism $F\to
KF/U$. Then we have
$$
F/D_{\widetilde{\mathfrak S^{n}}}(F)\subset (KF/U)^{*}.
$$
It means that the group $F/D_{\widetilde{\mathfrak S^{n}}}(F)$ is
nilpotent of the class nilpotency at most $(n-1)$. Let
$\gamma_n(F)$ is the $n$-th term of the lower central series of the
group $F$. So we have
$$
\gamma_n(F/D_{\widetilde{\mathfrak S^{n}}}(F))=1,
$$
and then
$$
\gamma_n(F)\subset D_{\widetilde{\mathfrak S^{n}}}(F).
$$


It is known  a fundamental result that for a free group $F$ there
is the equality $\gamma_{n}(F)=D_{n}(\mathbb Z,F)$ \cite{Magnus},
\cite{Grun}, \cite{Witt} (see also \cite{Roehl},
\cite{MikhPassi_book}). So for a field  $K$ of zero characteristic
we have
$$
\gamma_{n}(F)=D_{n}(K,F).
$$
Thus we have
$$
D_{n}(K,F)\subset D_{\widetilde{\mathfrak S^{n}}}(F).
$$

On the other hand,  $\mathfrak S^{n}\subset \widetilde{\mathfrak
S^{n}}$ and then $I_{\widetilde{\mathfrak S^{n}}}(KF)\subset
I_{\mathfrak S^{n}}(KF)$, where $I_{\widetilde{\mathfrak
S^{n}}}(KF)$ and $I_{\mathfrak S^{n}}(KF)$ are verbal ideals of
$KF$ for the varieties $\widetilde{\mathfrak S^{n}}$ and $\mathfrak
S^{n}$, respectively. Hence
$$
D_{\widetilde{\mathfrak S^{n}}}(F)\subset D_{\mathfrak
S^{n}}(F)=D_{n}(K,F).
$$

Thus
$$
D_{n}(K,F)=D_{\widetilde{\mathfrak S^{n}}}(F).
$$
Theorem is proved.
} 

We have the following
\begin{cor}
For a free group $F$  the dimension subgroup $D_{\mathfrak
S^{n}}(F)$ for the variety of representations of groups $\mathfrak
S^{n}$ coincides with the dimension subgroup
$D_{\widetilde{\mathfrak S^{n}}}(F)$ for the polynomial variety
$\widetilde{\mathfrak S^{n}}$.
\end{cor}

We can state

\begin{prob}\label{Prob_DimSub1}
Find a non-free group $G$ such that the dimension subgroup
$D_{n}(K,G)$ coincides with the dimension subgroup $D_{\mathfrak
X}(G)$ for some variety of representations of groups $\mathfrak
X\ne \mathfrak S^{n}$.
\end{prob}

\paragraph{5.3.2}

Let $\mathfrak X=\eta\Sigma$ be a polynomial variety of representations of groups.
Assume that the identities of the variety $\Sigma$ are known.

Let $G$ be a group, let $KG$ be the group algebra of  $G$
and let $I_{\mathfrak X}(KG)$ be the verbal ideal of identities
for the variety $\mathfrak X$ in the algebra $KG$. Note that we
have the equality $I_{\mathfrak X}(KG)=I_{\Sigma}(KG)$, where
$I_{\Sigma}(KG)$ is the verbal ideal of identities for the variety
of algebras $\Sigma$ in the algebra $KG$. The group $G$ admits the
representation $\rho=(KG/I_{\mathfrak X}(KG),G)$ from the variety
$\mathfrak X$ and the kernel of this representation is the dimension
subgroup of $G$ for the variety $\mathfrak X$, that is
$$
D_{\mathfrak X}(G)=\mbox{Ker}\rho=\{ g\in G\mid g-1\in
I_{\mathfrak X}(KG) \}.
$$
Moreover we have (see theorem \ref{Th_EqualDimSubgr}):
$$
D_{\mathfrak X}(G)=D_{\Sigma}(G).
$$

Along with the representation $\rho=(KG/I_{\mathfrak X}(KG),G)$ we
can consider the representation $\rho_1=(KG/I_{\mathfrak
X}(KG),KG)$. Since the representation $\rho$ belongs to the
variety $\mathfrak X=\eta\Sigma$ then the algebra
$\overline{KG}=KG/\mbox{Ker}\rho_1\cong \langle \overline G
\rangle $ is an algebra from the variety $\Sigma$, (where $\langle
\overline G \rangle$ is a linear span of the group $\overline
G=G/\mbox{Ker}\rho=G/D_{\Sigma}(G)$ in the algebra
$\mbox{End}(KG/I_{\mathfrak X}(KG))$ ). Thus,  $\overline
G\subset \overline{KG}$.

We have the following

\begin{prob}
We know the identities of the algebra $\overline{KG}\in \Sigma$.
What can we say about the identities of the group $\overline
G=G/D_{\Sigma}(G)$?
\end{prob}

We have solution of this problem for the variety $\Sigma$ of
algebras with the nilpotent derived ideal.

Let $\mathfrak N_{n-1}$ be the variety of nilpotent groups of
class $(n-1)$, and let $\mathfrak A$ be a variety of
abelian groups.

\begin{prop}
If $\Sigma$ is a variety of algebras with the nilpotent derived ideal
of class $n$ and a group $G$ admits a representation from the
variety $\eta\Sigma$ then the group $G/D_{\Sigma}(G)$ belongs to
the variety of groups $\mathfrak N_{n-1}\cdot \mathfrak A$.
\end{prop}

\dok{
 Let $\mathfrak X=\eta\Sigma$. In the
section~\ref{S4_VarOfReprGrVarAlg} it was showed that
$\eta\Sigma=\mathfrak S^{n}\times \mathfrak A=(\omega\mathfrak
A)^{n}$.

Let $G$ be a group and let $V=KG/I_{\Sigma}(KG)$ then there is the
representation $(V,G)$ from the variety $\mathfrak X$. Let
$\overline G=G/D_{\Sigma}(G)$. Since the representation
$(V,\overline G)$ belongs to the variety $\mathfrak X=\mathfrak
S^{n}\times \mathfrak A$ then there is a normal subgroup
$\overline H\triangleleft \overline G$ such that there is a
representation $(V,\overline H)$ from the variety $\mathfrak
S^{n}$ and the group $\overline G/\overline H$ is abelian.

The representation $(V,\overline H)$ is faithful (as a
subrepresentation of the faithful representation $(V,\overline
G)$) and it belongs to the variety $\mathfrak S^{n}$, so according
to the Kaloujnine's theorem \cite{Kaloujnine} the group $\overline
H$ is $(n-1)$-nilpotent. Moreover $\overline G/\overline H$ is
abelian and the group $\overline G$ belongs to the variety
$\mathfrak N_{n-1}\cdot \mathfrak A$.
} 

Let $\overline G=G/D_{\Sigma}(G)$ be a group from the variety of
groups $\Theta$. Let $B$ be the verbal subgroup of the group $G$ for
the variety $\Theta$. Then $G/B\in \Theta$ and $B\subset
D_{\Sigma}(G)$.

We can consider
\begin{prob}
What can we say about the structure of the group
$D_{\Sigma}(G)/B$?
\end{prob}

Let $B$ be the verbal  subgroup of the group $G$ for the variety
$\mathfrak N_{n-1}\cdot \mathfrak A$. Let $\Sigma$ be the variety of
algebras with the nilpotent derived ideal of the class $n$. We have the
following

\begin{prop}
Let a group $G$ admit a representation in the variety of
representations of groups $\eta\Sigma$. If the group $G/B$ is
torsion-free then $D_{\Sigma}(G)=B$.
\end{prop}

\dok{
 As it was noted above $B\subset D_{\Sigma}(G)$.
 We will show that $D_{\Sigma}(G)\subset B$.

 Let $G_1=G/B$. Since $B\subset D_{\Sigma}(G)$
 and $G/D_{\Sigma}(G)\in \mathfrak{N}_{n-1}\cdot \mathfrak A$
 then the group $G_1\in \mathfrak
N_{n-1}\cdot \mathfrak A$. So there is a normal subgroup
$G_2\triangleleft G_1$ such that $G_2$ is $(n-1)$-nilpotent and
the group $G_1/G_2$ is abelian. Since the group $G_1$ is
torsion-free then the subgroup $G_2$ is also torsion-free and
$(n-1)$-nilpotent. According to the result of
Maltsev~\cite{Malcev2} (see also \cite{PV})
if $K$ is a field of a characteristic zero then the class of all
groups admitting a faithful representation in the variety
$\mathfrak S^{n}$ coincides with the class $\mathfrak N_{n-1,0}$
of all torsion-free nilpotent groups. So  there is a
$K$-module $U$ such that the group $G_2$ admits a faithful
representation $(U,G_2)$ from the variety $\mathfrak S^{n}$.

So we have the group $G_1$ and the representation $(U,G_2)$ of the
subgroup $G_2$ of the group $G_1$. Here we can use induced representations.
Let $U_1=U\otimes_{KG_2} KG_1$, let
$(U_1,G_1)$ be a representation induced by the representation
$(U,G_2)$. So for every $u\in U_1$ and every $g\in G_1$ we have
$$
u_1\circ g=(u\otimes \sum_{i}\alpha_ig_i)\circ g= (u\circ
g)\otimes \sum_{i}\alpha_ig_i=u\otimes \sum_{i}\alpha_igg_i,
$$
where $u\in U$, $\alpha_i\in K$, $g_i\in G_1$. Since the
representation $(U,G_2)$ is faithful  the representation
$(U_1,G_1)$ is also faithful.
Let $h\in G_2$ then
$$
(U\otimes_{KG_2} KG_1)\circ h=(U\circ h)\otimes_{KG_2} KG_1.
$$
The group $G_2$ acts $n$-unitriangularly on $U$. Then $G_2$ acts
$n$-unitriangularly on $U\otimes_{KG_2} KG_1$ and the representation
$(U_1,G_2)$ belongs to the variety $\mathfrak S^{n}$.

Thus we have that the group $G_1$ has the normal subgroup $G_2$
such that $G_1/G_2$ is abelian and $(U_1,G_2)\in \mathfrak S^{n}$,
it means that the faithful representation $(U_1,G_1)$ belongs to
the variety $\mathfrak S^{n}\times \mathfrak A$.

With the representation $(U_1,G_1)=(U_1,G/B)$ we  consider the
representation $(U_1,G)$. It belongs to the variety $\mathfrak
S^{n}\times \mathfrak A$ and the kernel of this representation is
the subgroup $B$.
Since $D_{\Sigma}(G)=D_{\mathfrak X}(G)$ is the intersection of
kernels of all representations of the group $G$ from the variety
$\mathfrak X=\eta\Sigma$ then $D_{\Sigma}(G)\subset B$.

So we have proved that $D_{\Sigma}(G)=B$.
 } 

Remind that the variety $(\omega\mathfrak A)^{n}$ is generated by
the representation $(K^{n}, T_{n}(K))$ and the variety $\mathfrak
S^{n}$ is generated by the representation $(K^{n}, UT_{n}(K))$. So
the classical dimension subgroup is related to  the group of
unitriangular matrices $UT_{n}(K)$ and dimension subgroup for the
polynomial variety $(\omega\mathfrak A)^{n}$ is related to the
group of triangular matrices $T_{n}(K)$.




\begin{thebibliography}{99}

\bibitem{AladovaPlotkin}
E.V. Aladova, B.I. Plotkin, $PI$-groups and $PI$-representations
of groups. Fundamental and applied mathematics, 2009, to appear.
Arxiv: math RA: 0905.3651.


\bibitem{Bergman-Lewin}
G.M.~Bergman, J.~Lewin, The semigroup of ideals of a fir is
(usually) free, J. London Math. Soc., {\bf 11},(1975), no.~1,
p.~21--31.

\bibitem{Birkhoff}
G.~Birkhoff, On the structure of abstract algebras, Proc.
Cambridge Philos. Soc., {\bf 31}, (1935), p.~433--454.

\bibitem{Fox}
R.~Fox, Free differential calculus. I. Derivation in the free
group ring. Ann. of Math., {\bf 57}, (1953), no.~2, p.~547--560.


\bibitem{Grinberg}
A.C.~Grinberg, Varieties of stable linear representations.
(Russian). Latvian math. yearbook. Izdat. ''Zinatne'', Riga, {\bf
9}, (1971), p.~39--46.

\bibitem{Grun}
 O.~Gr\" un, \" Uber die Faktogruppen freier Gruppen I.
Deutsche Math. (Jahrgang 1), {\bf 6}, (1936), p.~772--782.

\bibitem{Gupta}
N.~Gupta, The dimension subgroup conjecture. Bull. London Math.
Soc., {\bf 22}, (1990), no.~5, p. 453--456.




\bibitem{Kaljulaid}
U.E.~Kaljulaid, Triangular pruducts of representations of
semigroups and associtive algebras, Uspehi Math. Nauk, {\bf 32},
no. 4, (1977), p.~ 253--254.


\bibitem{Kaloujnine}
L.~Kaloujnine, \" Uber gewisse Beziehungen zwischen einer Gruppe
und ihren Automorphismen. (German), Beriliner Math. Tagung,
(1953), pp. 164--172.





\bibitem{Krop}
L.E.~Krop, Solvable varieties of pairs. (Russian) Latvian. Mat.
yearbook.  {\bf 18}, (1976), p.~64--80.



\bibitem{MacLane}
S.~Mac Lane, Categories for the working mathematican, Second
edition. Graduate text in mathematics, 5. Springer-Verlag, New
York, (1998), xii+314~pp.

\bibitem{Magnus}
W.~Magnus, \" Uber Beziehungen zwieschen h\" oheren Kommutatoren.
J. reine angew. Math., {\bf 177}, (1937), p.~105--115.



\bibitem{Malcev}
A.I.~Malcev, Algebraic systems. Posthumous edition, edited by D.
Smirnov and M. Taiclin. Translated from the Russian by B. D.
Seckler and A. P. Doohovskoy. Die Grundlehren der mathematischen
Wissenschaften, Band 192. Springer-Verlag, New York-Heidelberg,
1973.

\bibitem{Malcev2}
A.I.~Malcev, Generalized nilpotent algebras and their associated
groups. (Russian) Mat. Sbornik N.S. 25(67), (1949). 347--366.



\bibitem{MikhPassi}
R.~Mikhailov, I.B.S.~Passi, Augmentation powers and group
homology. J. Pure Appl. Algebra {\bf 192} (2004), no. 1-3,
p.~225--238.

\bibitem{MikhPassi_book}
R.~Mikhailov, I.B.S.~Passi, Lower central and dimension series of
groups. Lecture Notes in Mathematics, 1952. Springer-Verlag,
Berlin, 2009. xxii+346 pp.

\bibitem{MikhPassi_D4}
R.~Mikhailov, I.B.S.~Passi, The quasi-variety of groups with
trivial fourth dimension subgroup. J. Group Theory, {\bf 9},
(2006), p.~369--381.


\bibitem{Neumann_3}
B.H.~Neumann, H.~Neumann, P.M.~Neumann, Wreath products and
varieties of groups, Math. Ztschr., {\bf 80}, (1962), no.~1,
p.~44--62.


\bibitem{Passi}
I.B.S.~Passi, Group rings and their augmentation ideals. Lecture
Notes in Mathematics, {\bf 715}. Springer, Berlin, (1979). vi+137
pp.

\bibitem{Plotkin_Autom}
B.I.~Plotkin, Group of automorphisms of algebraic system.
Translated from the Russian by K.A.~Hirsh. Wolters-Noordhoff
Publishing, Groningen, (1972), xviii+502~pp.


\bibitem{Plotkin_Notes}
B.I.~Plotkin, Notes on Engel groups and Engel elements in groups.
Some generalizations, (English, Russian summary), Izv. Ural. Gos.
Univ. Mat. Mekh., {\bf 36}, (2005), no.~7, p.~153--166,
p.~192--193.

\bibitem{Plotkin_VarAndVarOfPairs}
B.I.~Plotkin, Varieties of groups, and varieties of pairs that are
connected with representations of groups. (Russian), Sibirsk. Mat.
Z. {\bf 13}, (1972), no. 5, p.~1030--1053.


\bibitem{Plotkin_VarAndOuasivar}
B.I.~Plotkin,  The varieties and quasivarieties that are connected
with group representations, (Russian), Dokl. Akad. Nauk SSSR, {\bf
196} (1971), no.~3, p.~527--530.

\bibitem{Plotkin-Grinberg}
B.I.~Plotkin, A.S.~Grinberg, Semigroups of varieties that are
connected with group representations. (Russian), Sibirsk. Mat. Z.,
{\bf 13}, (1972), no. 5, p. 841--858.


\bibitem{PV}
B.I.~Plotkin, S.M.~Vovsi, Varieties of group representations.
General theory, connections and applications, (Rusian),
''Zinatne'', Riga, (1983), 339~pp.







\bibitem{Rips}
E.~Rips, On the fourth integer dimension subgroup. Israel J.
Math., {\bf 12}, (1972), p.~342--346.

\bibitem{Roehl}
F.~R\" ohl, Review and some critical comments on a paper of Grun
concerning the dimension subgroup conjecture. Bol. Soc. Brasil.
Mat. {\bf 16}, no. 2, (1985), p.~11--27.





\bibitem{Shmel'kin}
A.L.~Smel'kin, Semigroup of group varieties, (Rusian), Dokl. Akad.
Nauk SSSR, {\bf 149}, (1963), no.~3, p.~543--545.

\bibitem{Vovsi}

S.M.~Vovsi, Triangular products of group representations and their
applications. Progress in Mathematics, 17. Birkhauser, Boston,
Mass., 1981. ix+127 pp.

\bibitem{Witt}
E.~Witt, Treue Darstellung Lieschen Ringe. J. Reine angew. Math.,
{\bf 177}, (1937), p.~152--160.




\end{thebibliography}
\end{document}